\documentclass[12pt]{article}
 
\usepackage{bbm, amsmath, amsfonts, amscd, latexsym, amsthm, amssymb, graphicx, euscript, mathrsfs}
\usepackage[nottoc,notlot,notlof]{tocbibind}
\usepackage{tikz}
\usetikzlibrary{patterns}
\providecommand{\keywords}[1]{\textbf{Key words:} #1}
\newtheorem{theor}{\hspace{1cm}{\sc Theorem}}[section]
\newtheorem{utver}[theor]{\hspace{1cm}{\sc Proposition}}
\newtheorem{sledst}[theor]{\hspace{1cm}{\sc Corollary}}
\newtheorem{lemma}[theor]{\hspace{1cm}{\sc Lemma}}

\newtheorem*{utver*}{\hspace{1cm}{\sc Proposition}}
\theoremstyle{definition}
\newtheorem{defin}[theor]{\hspace{1cm}{\sc Definition}}
\newtheorem{exa}[theor]{\hspace{1cm}{\sc Example}}
\newtheorem{rem}[theor]{\hspace{1cm}{\sc Remark}}

\newcommand{\Vol}{\mathop{\rm Vol}\nolimits}

\newcommand{\determ}{\mathop{\rm {det}}\nolimits}
\newcommand{\codim}{\mathop{\rm codim}\nolimits}

\newcommand{\rk}{\mathop{\rm rk}\nolimits}

\newcommand{\conv}{\mathop{\rm conv}\nolimits}

\newcommand{\MV}{\mathop{\rm MV}\nolimits}

\newcommand{\supp}{\mathop{\rm supp}\nolimits}
\newcommand{\val}{\mathop{\rm val}\nolimits}                                                                                             
\newcommand{\tzero}{{\mbox{\bfseries \large 0}}}
\newcommand{\tid}{{\mbox{\bfseries\large 1}}}
\newcommand{\Mod}[1]{\ (\text{mod}\ #1)}
\def\R{\mathbb R}

\def\Z{\mathbb Z}
\def\F{\mathbb F}

\def\C{\mathbb C}
\def\CC{({\mathbb C}\setminus 0)}

\def\T{\mathbb T}

\def\K{\mathbb K}
\def\newton{\EuScript{N}}

\def\resul{\mathscr{R}}

\newcommand*\samethanks[1][\value{footnote}]{\footnotemark[#1]}

\begin{document}

\title {Signs of the Leading Coefficients of the Resultant}
\author{Arina Arkhipova\thanks{National Research University Higher School of Economics \newline The article was prepared within the framework of the Academic Fund Program at the National Research University Higher School of Economics (HSE) in 2016-2017(grant N 16-01-0069) and supported within the framework of a subsidy granted to the HSE by the Government of the Russian Federation for the implementation of the Global Competitiveness Program.}, Alexander Esterov\samethanks}
\date{}

\maketitle{}

\begin{abstract}

We construct a certain $\F_2$-valued analogue of the mixed volume of lattice polytopes. 
This 2-mixed volume cannot be defined as a polarization of any kind of an additive measure, 
or characterized by any kind of its monotonicity properties, because neither of the two makes sense over $\F_2$. In this sense, the convex-geometric nature of the 2-mixed volume remains unclear.

On the other hand, the 2-mixed volume seems to be 
no less natural and useful than the classical mixed volume -- in particular, it also plays an important role in algebraic geometry. As an illustration of this role, we obtain a closed-form expression in terms of the 2-mixed volume to compute the signs of the leading coefficients of the resultant, which were by now explicitly computed only for some special cases.
\end{abstract}

\\

{\bf \Large Introduction}\\

The classical resultant was initially studied by Sylvester (1853), and later extended to the case of a system of $n$ homogeneous polynomials in $n$ 
variables by Cayley (1948) and Macaulay (1902).
In the 1990s, the advances in several fields, such as symbolic algebra and multivariate hypergeometric functions, revived the interest in resultants. Sparse resultants were introduced and studied by Agrachev, Gelfand, Kapranov, Zelevinsky, and Sturmfels (see e.g. \cite{gkz}). In particular, in \cite{sturmfels}, Sturmfels gives an explicit combinatorial construction of the Newton polytope of the sparse resultant, and proves that 
the leading coefficient of the resultant with respect to an arbitrary monomial order is equal to $\pm1$. However, the signs of such coefficients have been computed explicitly only for some special cases so far, although the general answer might be useful for the purposes of real algebraic geometry. 

In our work, we construct the {\it 2-mixed volume} (Definition \ref{defmv2}), which is an analogue of the classical mixed volume of convex lattice polytopes taking values in $\F_2$. Besides that, we express the signs of the leading coefficients of the sparse resultant in terms of the 2-mixed volume of certain tuples of polytopes (Theorem \ref{restheo}). 

The 2-mixed volume is a symmetric and multilinear function of lattice polytopes (Proposition \ref{propmv2}). However, its convex-geometric nature remains unclear, because we cannot define it as a polarization of any kind of an additive measure, or characterize it by any kind of its monotonicity properties.

Our explicit formula for the 2-mixed volume employs the so-called 2-determinant, that is, the unique nonzero multilinear function of $n+1$ vectors in the $n$-dimensional vector space over the field $\F_2$ which ranges in $\F_2$, remains invariant under all linear transformations, and equals zero whenever the rank of the $n+1$ vectors is less than $n$.
This function implicitly appeared in the context of the class field theory for multidimensional local fields by Parshin and Kato (see e.g. Remark 1 in Section 3.1 of \cite{parshin}, which is probably the first occurence of the 2-determinant in the literature). Later this notion was explicitly introduced in full generality by A.Khovanskii in \cite{det2} for the purpose of his multidimensional version of the Vieta formula (i.e. the computation of the product in the group $\CC^n$ of all the roots for a system of $n$ polynomial
equations with sufficiently generic Newton polytopes).

The algebro-geometric part of our work is an extension of related results by A. Khovanskii. In particular, our notion of the 2-mixed volume is the result of our effort to provide an invariant interpretation of the sign in Khovanskii's multivariate version of the Vieta formula, and to relax the genericity assumptions on the Newton polytopes in this formula. 

The convex-geometric part of our work employs the techniques of tropical geometry to prove the existence of the 2-mixed volume (Theorem \ref{main}).\\

{\bf\Large Structure of the paper}\\

The paper is organized as follows. In Section 1, we recall some necessary facts and notation concerning convex and tropical geometry. 

Section 2 is devoted to the notion of the {\it 2-mixed volume}. First, we recall the definition and the basic properties of the {\it 2-determinant}, and use it to define the so-called {\it 2-intersection number} of tropical hypersurfaces. Then we show that the 2-intersection number depends only on the Newton polytopes of the hypersurfaces, which yields a well-defined function of lattice polytopes --- the so-called  {\it 2-mixed volume}. 

Section 3 concerns the multivariate Vieta's formula which expresses the product of roots for a polynomial system of equations in terms of the 2-mixed volume of its Newton polytopes. 

In Section 4, we first recall the definition of the sparse mixed resultant, then we compute the signs of the leading coefficients of the resultant reducing this problem to finding the product of roots for a certain system of equations (see Subsection \ref{SignsRes}, Theorem \ref{restheo}). \\

{\small
\tableofcontents}
\section{Preliminaries}
\subsection{Some Definitions, Notation and Basic Facts}
Here we introduce some basic notation, definitions and facts that will be used throughout this paper. For more details, we refer the reader to the works \cite{bernstein}, \cite{EsKhov} and \cite{tropsturmfels}.
\subsubsection{Laurent Polynomials and Newton Polytopes}
\begin{defin}
A {\it Laurent polynomial} $f(x_1,\ldots, x_n)$ in the variables $x_1,\ldots, x_n$ over a field $\F$ is a formal expression $$f(x_1,\ldots, x_n)=\sum_{(a_1,\ldots, a_n)\in\Z^n} c_{a_1,\ldots, a_n} x_1^{a_1}\ldots x_n^{a_n},$$ where the coefficients $c_{a_1,\ldots, a_n}$ belong to $\F$ and only finitely many of them are non-zero. We denote the ring of Laurent polynomials in $n$ variables with coefficients in $\F$ by $\F[x_1^{\pm 1},\ldots, x_n^{\pm 1}]$. 
\end{defin}

\begin{rem}
Throughout this paper, we use the multi-index notation, i.e., for $a=(a_1,\ldots, a_n)\in \Z^n$, instead of $c_{a_1,\ldots, a_n} x_1^{a_1}\ldots x_n^{a_n}$, we will use the expression $c_a x^a$. 
\end{rem}

\begin{defin}
Let $f(x)=\sum_{a\in \Z^n} c_a x^a$ be a Laurent polynomial. The {\it support of $f$} is the set $\supp(f)\subset \Z^n$ which consists of all points $a\in \Z^n$ such that the corresponding coefficient $c_a$ of the polynomial $f$ is non-zero. 
\end{defin}

\begin{defin}
The {\it Newton polytope of $f(x)$} is the convex hull of $\supp(f)$ in $\R^n$, i.e., the minimal convex lattice polytope in $\R^n$ containing the set $\supp(f)$. We denote the Newton polytope of a polynomial $f(x)$ by $\newton(f)$. 
\end{defin}

\begin{defin} \label{minkowski}
For a pair of subsets $A, B\subset\R^n$, their {\it Minkowski sum} is defined to be the set $A+B=\{a+b\mid a\in A, b\in B\}$. 
\end{defin}

The following important fact provides a connection between the operations of the Minkowski addition and multiplication in the ring of Laurent polynomials.

\begin{utver} 
For a pair of Laurent polynomials $f,g$, we have the following equality: $$\newton(fg)=\newton(f)+\newton(g).$$
\end{utver}

\begin{defin}\label{suppface}
Let $A\subset\R^m$ be a convex lattice polytope and $\ell\in(\R^*)^m$ be a covector. Consider $\ell$ as a linear function, and denote by $\ell\mid_A$ its restriction to the polytope $A$. The function $\ell\mid_A$ attains its maximum at some face $\Gamma\subset A$. This face is called {\it the support face} of the covector $\ell$ and is denoted by $A^{\ell}$.
\end{defin}

\subsubsection{Cones and Fans}
\begin{defin}
A {\it polyhedral cone} $C$ in $\R^n$ is the positive hull of a finite subset $S=(s_1,\ldots, s_m)\subset\R^n$, i.e., $C=\{\sum_{i=1}^m \lambda_i s_i\mid \lambda_i\geqslant 0\}$. 
\end{defin}

\begin{defin}
Let $C\subset\R^n$ be a polyhedral cone. A face of $C$ is the intersection $\{l=0\}\cap C$ for a linear function $l:\R^n\to\R$ such that $C\subset \{l\geqslant0\}$. 
\end{defin}

\begin{defin}
A {\it polyhedral fan} $\Sigma$ in $\R^n$ is a collection of polyhedral cones which satisfies the following properties: every face of every cone from $\Sigma$ is an element of $\Sigma$, and for any pair of cones $S_1, S_2\in\Sigma$, $S_1\cap S_2$ is a face for both $S_1$ and $S_2$.
\end{defin}

Throughout this paper, will mostly deal with the following special case of a polyhedral fan.

\begin{defin}
For a convex polytope $A\subset R^m$ and a  face $\Gamma\subset A$, we define the {\it normal cone} $N(\Gamma)$ to be the union of all the covectors $\ell\in (\R^*)^m$ such that the support face $A^{\ell}\subset A$ (see Definition \ref{suppface}) contains $\Gamma$.

Then, the {\it normal fan} of the polytope $A$ is the collection $N(A)=\{N(\Gamma)\,|\,\Gamma\subset A\}$ over all the faces $\Gamma\subset A$. 
\end{defin}

\begin{defin}
Let $\Sigma$ be a polyhedral fan in $\R^n$. The {\it support}~$\supp(\Sigma)$ of $\Sigma$ is the union of all of its cones.
\end{defin}

\begin{defin}
Let $\Sigma_1, \Sigma_2$ be polyhedral fans. The {\it common refinement} $\Sigma_1\wedge \Sigma_2$ is defined to be the fan consisting of all the intersections $C_1\cap C_2$, where $C_i\in \Sigma_i$. 
\end{defin}

\begin{defin}
Let $P=(P_1,\ldots, P_m)$ be a tuple of polytopes in $\R^n$. A fan $\Sigma$ is said to be {\it compatible} with the tuple $P$, if each of its cones is contained in a cone of the common refinement~${N(P_1)\wedge N(P_2)\wedge\ldots\wedge N(P_m)}$ of the normal fans of the polytopes in $P$. 
\end{defin}

The following statement relates the Minkowski sums (see Definition \ref{minkowski}) to normal fans. 

\begin{utver}
 Let $P, Q\subset\R^n$ be polytopes. Then the following equality holds: $$N(P)\wedge N(Q)=N(P+Q).$$
\end{utver}

\subsubsection{The Mixed Volume and the Bernstein--Kushnirenko Formula}

\begin{defin}
Let $\gamma\neq 0$ in $(\R^*)^n$ be a covector and $f(x)$ be a Laurent polynomial with the Newton polytope $\newton(f)$. The {\it truncation} of $f(x)$ with respect to $\gamma$ is the polynomial $f^{\gamma}(x)$ that can be obtained from $f(x)$ by omitting the sum of monomials which are not contained in the support face ${\newton(f)^{\gamma}}$. 
\end{defin}

It is easy to show that for a system of equations $\{f_1(x)=\ldots=f_n(x)=0\}$ and an arbitrary covector $\gamma\neq 0$, the  system $\{f_1^{\gamma}(x)=\ldots=f_n^{\gamma}(x)=0\}$ by a monomial change of variables can be reduced to a system in $n-1$ variables at most. Therefore, for the systems with coefficients in general position, the ``truncated" systems are inconsistent in $\CC^n$. 

\begin{defin}
Let $(f_1,\ldots, f_n)$ be a tuple of Laurent polynomials. In the same notation as above, the system ${f_1(x)=\ldots=f_n(x)=0}$ is called {\it degenerate at infinity}, if there exists a covector $\gamma\neq 0$ such that the system $\{f_1^{\gamma}(x)=\ldots=f_n^{\gamma}(x)=0\}$ is consistent in $\CC^n$.
\end{defin}

\begin{defin} \label{mv}
Let $\mathscr P$ be the semigroup of all convex polytopes in $\R^n$ with respect to the Minkowski addition (see Definition \ref{minkowski}). The {\it mixed volume} is a unique function $$\MV\colon\underbrace{\mathscr{P}\times\ldots\times\mathscr{P}}_{\mbox{$n$ times}}\to\R$$ which symmetric, multilinear (with respect to the Minkowski addition) and which satisfies the following property: the equality $MV(P,\ldots, P)=\Vol(P)$ holds for every polytope~${P\in \mathscr P}$.
\end{defin}

The following theorem allows to compute the number of roots for a non-degenerate polynomial system of equations in terms of the mixed volume of its polynomials. 

\begin{theor}[Bernstein-Kushnirenko formula, \cite{bernstein}]
The number of roots for a polynomial system of equations $\{f_1(x)=\ldots=f_n(x)=0\}$ in $\CC^n$ that is non-degenerate at infinity counted with multiplicities is equal to $n!\MV(\newton(f_1),\ldots, \newton(f_n))$.
\end{theor}

\subsection{A Very Little Bit of Tropical Geometry}
\begin{defin}
We define the {\it tropical semifield} $\T=\R \cup \{-\infty\}$ to be the set of real numbers with $-\infty$ equipped with the following arithmetic operations:
$$
\alpha\oplus \beta=\begin{cases}
\max(\alpha,\beta),&\text{if $\alpha\neq \beta$;}\\
[-\infty,\alpha],&\text{if $\alpha=\beta$;}
\end{cases}
$$
$$
\alpha\odot \beta = \alpha+\beta.
$$
\end{defin}
\begin{rem}
Tropical addition and multiplication have the identity elements: $\tzero=-\infty$ and $\tid =0$.
\end{rem}
\begin{rem}
It is easy to check the following facts:
\begin{itemize}
\item $\T$ is a commutative multivalued semigroup (a hypersemigroup) with respect to addition;
\item $\T\setminus\{-\infty\}$ is a commutative group with respect to multiplication;
\item in $\T$ we have the distribution law: $\forall \alpha,\beta,\gamma \in \T~ \alpha\odot (\beta\oplus \gamma) = \alpha\odot\beta\oplus\alpha\odot\gamma$.
\end{itemize}
\end{rem}

Having defined tropical arithmetic operations, we can consider tropical polynomials.
\begin{defin}\label{troppol}
Let $A \subset \Z^n$ be finite and ${\forall a \in A} ~c_a \in \T$. Then a {\it tropical polynomial} is given by
$$
f(x) =\bigoplus_{a\in A} c_a\odot x^{\odot a},
$$ where $x\in \T^n$.
\end{defin}
\begin{defin}
In the notation of \ref{troppol}, the {\it support} of the polynomial $f(x)$ is the set $\supp(f)=\{a\in A\mid c_a\neq \tzero\}$ (see ~\ref{troppol}). The {\it Newton polytope} $\newton(f)$ of $f(x)$ is defined as the convex hull of $\supp(f)$ in $\R^n$. We denote by $\vert supp(f)\vert$ and $|\newton(f)|$ the cardinality of the set $\supp(f)$ and $\newton(f)\cap\Z^n$, respectively.
\end{defin}
\begin{defin}
We denote by $=_{\mathbb T}$ {\it the multivalued equality sign}: $f(x_0)=_{\mathbb T}\tzero$ for some $x_0\in \T^n$ means that~${\tzero=-\infty}$ belongs to the image $f(x_0)$.
\end{defin}
\begin{defin}\label{weight1}
Let $f(x)=\bigoplus_{a\in A} c_a\odot x^{\odot a}$ be a tropical polynomial. Consider the set ${H=\{x_0\in \T^n\mid f(x_0) =_{\T}\tzero\}}$ (or, equivalently, the set of points $x_0\in \T^n$ such that there exist $a_1\neq a_2 \in A$ satisfying $c_{a_1}+x_0\cdot a_1=c_{a_2}+x_0\cdot a_2=\max f(x_0)$).

A point ${x_0\in~H}$ is said to be {\it smooth} if the points $a_1$ and $a_2$ above are uniquely defined (up to the transposition). 
In this case, the {\it weight} of $H$ at $x_0$ is defined as the integer length of the vector $a_1-a_2$ (i.e. the g.c.d. of its coordinates).

\end{defin}
\begin{defin}\label{trophyp1}
The {\it tropical hypersurface} defined by $f$ is the set $H$ whose facets (i.e., the connected components of the smooth part of $H$) are equipped with weights (see~\ref{weight1}).
\end{defin}

As a set, the tropical hypersurface $f =_{\T}\tzero$ is just the set of points where $f$ is not smooth, i. e., the corner locus of this convex piecewise linear function.
 
See, for example, \cite{trop1} and \cite{tropsturmfels} for a more detailed introduction into tropical geometry.

\begin{defin} The local support set of a polynomial $f(x)=\bigoplus_{a\in A} c_a\odot x^{\odot a}$ at a point $x_0$ is the set $\supp_{x_0}(f)$ of all $a$ such that $c_{a}+x_0\cdot a=\max f(x_0)$. The local Newton polytope $\newton_{x_0}(f)$ is the convex hull of the local support set $\supp_{x_0}(f)$.
\end{defin}
Note that $x_0$ belongs to the hypersurface $H=\{f=_{\T}\tzero\}$ if and only if $\supp_{x_0}(f)$ consists of more than one point, and is smooth if $\supp_{x_0}(f)$ consists of two points. Also note that the local Newton polytope (in contrast to the local support set) depends only on the hypersurface $H$, and not on its defining equations $f$, so we shall also denote it by $\newton_{x_0}(H)$. 

\section{2-mixed volume}\label{2-volume}
This Section is devoted to the notion of the {\it 2-mixed volume} of lattice polytopes. In Subsection \ref{defdet2}, the definition and some basic properties of the 2-determinant are provided. Subsection \ref{2ind} concerns the definition of the {\it 2-intersection number} for a tuple of tropical hypersurfaces. In Subsections \ref{ideamain}, \ref{buildwall} and \ref{passwall} we prove Theorem \ref{main} stating that under some assumptions on the Newton polytopes of the tropical hypersurfaces, the 2-intersection number depends not on the hypersurfaces, but on their Newton polytopes, which yields a well-defined function of lattice polytopes that takes values in $\mathbb F_2$ -- the {\it 2-mixed volume}.
\subsection{Analog of the determinant for $n+1$ vectors in an $n$-dimensional space over the field $\mathbb F_2$}\label{defdet2}

\begin{defin}
We define $\det_2$ to be the function of $n+1$ vectors in an $n$-dimensional linear space over $\mathbb F_2$, that takes values in $\mathbb F_2$ and satisfies the following properties:
\begin{itemize}
\item $\det_2(k_1,\ldots,k_{n+1})$ is equal to zero, if the rank of the collection of vectors $k_1,\ldots,k_{n+1}$ is smaller than $n$;
\item $\det_2(k_1,\ldots,k_{n+1})$ is equal to $\lambda^1+\ldots+\lambda^{n+1}+1$, if the vectors $k_1,\ldots,k_{n+1}$ are related by the unique relation $\lambda^1k_1+\ldots+\lambda^{n+1}k_{n+1}=0$.
\end{itemize}
\end{defin}
\begin{lemma}\label{det2}
The function $\det_2$
\begin{enumerate}
\item is $\mathrm GL_n(\mathbb F_2)$-invariant, i.e. for any linear transformation $A\in\mathrm GL_n(\mathbb F_2)$ the equality $\det_2(k_1,\ldots,k_{n+1})=\det_2(Ak_1,\ldots,Ak_{n+1})$ holds;
\item is multilinear.
\end{enumerate}
\end{lemma}
\begin{theor}\cite{det2}
There exists a unique nonzero function $\det_2$ which satisfies the properties of g ~\ref{det2}.
\end{theor}
\begin{theor}\cite{det2}\label{formuladet2}
In coordinates the function $\det_2$ can be expressed by the formula
\begin{equation*}
{\determ_2}(k_1,\ldots,k_{n+1})=\sum_{j>i}\Delta_{ij},
\end{equation*}
where $\Delta_{ij}$ is the determinant of the $n \times n$ matrix whose first $n-1$ columns represent the sequence of vectors $k_1,\ldots k_{n+1}$ from which the vectors with the indices $i$ and $j$ are deleted, and the last column is the coordinate-wise product of the vectors $k_i$ and $k_j$.
\end{theor}
Let $k_1,\ldots, k_{n+1}$ be a tuple of vectors such that $\rk(k_1,\ldots, k_{m+1})=m$ for some $m<n$. Then, there exists a natural projection $\pi\colon\Z^n\to\Z^n/\langle k_1,\ldots, k_{m+1}\rangle$. The statement below easily follows from Theorem \ref{formuladet2} and the well-known formula for computing the upper triangular block matrix determinant.
\begin{sledst}[upper triangular block matrix 2-determinant]\label{blockdet2}
In the same notation as above, the following equality holds: 
\begin{equation*}
\determ_2(k_1,\ldots, k_{n+1})=\determ_2(k_1,\ldots,k_{m+1})\cdot\det(\pi(k_{m+2}),\ldots,\pi(k_{n+1})).
\end{equation*}
\end{sledst}

\subsection{2-intersection number}\label{2ind}
\begin{defin}
Let $H_1, \ldots, H_n$ be tropical hypersurfaces. We say that $H_1, \ldots, H_n$ intersect {\it transversely} (denote by $H_1\pitchfork \ldots\pitchfork H_n$), if $|H_1\cap H_2\cap\ldots \cap H_n|<\infty$ and all the points~$x\in H_1\cap H_2\cap \ldots \cap H_n$ are smooth for every $H_i$ (see Definition \ref{weight1}).
\end{defin}
\begin{defin}\label{1-index}
Let $H_1\pitchfork \ldots\pitchfork H_n$ be a transverse tuple of tropical hypersurfaces. The {\it intersection number} $\iota(H_1,\ldots, H_n)\in \Z$ is the sum
\begin{equation}\label{eq2ind}
\iota(H_1,\ldots, H_n)\stackrel{\mathrm{def}}{=}\sum\limits_{x\in H_1\cap H_2\cap \ldots \cap H_n} {\mathrm{det}}(\newton_x(H_1),\ldots \newton_x(H_n)).
\end{equation}
\end{defin}
It is well known that the intersection number of tropical hypersurfaces depends only on their Newton polytopes (and coincides with the mixed volume of the Newton polytopes). This fact is often referred to as the tropical Bernstein--Kushnirenko formula. We shall need the following $\F_2$-verison of the intersection number.
\begin{defin}\label{2-index}
Consider an arbirtary point $\zeta\in~\mathbb Z^n$. Let $H_1\pitchfork \ldots\pitchfork H_n$ be a transverse tuple of tropical hypersurfaces. We define the {\it 2-intersection number} $\iota_2(H_1,\ldots, H_n;\zeta)\in \F_2$ as follows:
\begin{equation}\label{eq2ind}
\iota_2(H_1,\ldots, H_n;\zeta)\stackrel{\mathrm{def}}{=}\sum\limits_{x\in H_1\cap H_2\cap \ldots \cap H_n} {\mathrm{det}_2}(\newton_x(H_1),\ldots \newton_x(H_n),\zeta).
\end{equation}
\end{defin}
Unfortunately, in general, the 2-intersection number does depend on the tropical hypersurfaces, and not only on their Newton polytopes. However, this dependence disappears if the Newton polytopes themselves are in general position in a sense that we describe below.
\begin{defin}
Let $P\subset \mathbb R^n$ be a polytope or a finite set. We define the {\it support face} of a covector $v\in (\R^n)^*$ to be the maximal subset of $P$ on which $v\mid_P$ attains its maximum. We shall denote this face by $P^v$.
\end{defin}

\begin{defin}
A finite set $P\subset \Z^n$ is called a {\it $2$-vertex}, if for any pair of points $p_1, p_2\in P,~p_1\equiv p_2 \Mod{2}$ (i.e., the corresponding coordinates of the points $p_1, p_2$ are of the same parity). A lattice polytope is called a {\it $2$-vertex}, if the set of its vertices is  a {\it $2$-vertex}.
\end{defin}

\begin{defin}\label{devwrt}
Let $P_1, \ldots, P_n$ be convex lattice polytopes in $\R^n$ or finite sets in $\Z^n$, and $\zeta$ be a point in $\Z^n$. The tuple $P_1,\ldots, P_n$ is said to be {\it 2-developed with respect to} $\zeta$ if, for any covector $v\in(\Z^n)^*$ such that $v(\zeta)\not\equiv 0 \mod 2$, there exists $i\in \{1,\ldots, n\}$ such that the support face $P_i^v$ is a $2$-vertex.
\end{defin}

\begin{defin}\label{prickly}
A tuple $P=(P_1,\ldots, P_n)$ of convex lattice polytopes is said to be {\it $\zeta$--prickly}, if for any covector $v\in (\R^*)^n$ such that $v(\zeta)\neq 0$, there exists $i\in \{1,\ldots, n\}$ such that the support face $P_i^{v}$ is a vertex. 
\end{defin}

\begin{rem}
Obviously, if a tuple $P$ is $\zeta$--prickly, then it is $2$--developed with respect to $\zeta$.
\end{rem}

\begin{theor}\label{main}
Consider a point $\zeta\in\Z^n$ and finite lattice sets $P_1, \ldots, P_n$. Suppose that $P_1, \ldots, P_n$ are  2-developed with respect to $\zeta$. Then for any two tuples $(H_1,\ldots, H_n)$ and $(H'_1,\ldots, H'_n)$ of tropical hypersurfaces, whose equations are supported at $P_1, \ldots, P_n$, the 2-intersection numbers 
$\iota_2(H_1,\ldots, H_n;\zeta)$ and $\iota_2(H'_1,\ldots, H'_n;\zeta)$ coincide. 
\end{theor}

\begin{defin}\label{defmv2}
For a tuple of polytopes $P_1,\ldots,P_n$, 2-developed with respect to $\zeta\in\Z^n$, consider generic tropical hypersurfaces $H_1,\ldots,H_n$, such that the equation of $H_i$ is supported at the set of vertices of $P_i$. Then the function ${\MV_2\colon(P_1,\ldots,P_n;\zeta)\mapsto\iota_2(H_1,\ldots, H_n;\zeta)}$ is well-defined. We call it the {\it 2-mixed volume}.
\end{defin}
 
 \begin{utver} \label{propmv2}
 The function $\MV_2$ is symmetric and multiplinear with respect to the Minkowski summation of the arguments.
 \end{utver}
 
 $\vartriangleleft$ The symmetry is obvious. In order to prove the additivity $\MV_2(P,P_2,\ldots,P_n;\zeta)+\MV_2(Q,P_2,\ldots,P_n;\zeta)=\MV_2(P+Q,P_2,\ldots,P_n;\zeta)$ whenever the summands make sense, chose generic tropical polynomials $p,q,p_2,\ldots,p_n$ with the Newton polytopes $P,Q,P_2,\ldots,P_n$ respectively. Then the 2-intersection numbers in the tautological equality $\iota_2(p=0,p_2=0,\ldots,p_n=0;\zeta)+\iota_2(q=0,p_2=0,\ldots,p_n=0;\zeta)=\iota_2(p\cdot q=0,p_2=0,\ldots,p_n=0;\zeta)$ make sense and equal the corresponding 2-mixed volumes.
 $\vartriangleright$
 
\subsection{The Idea of the Proof of Theorem \ref{main}}\label{ideamain}
A tropical polynomial $\varphi$ corresponds to a point in $\T^{\vert\supp(\varphi)\vert}$. Namely, to every tropical polynomial we associate the collection of its coefficients. Therefore, tropical hypersurfaces defined by tropical polynomials with some fixed support $A\subset\Z^n$ can be considered as points in $\T^{\vert A\vert}$. 

Thus, given a tuple $(A_1,\ldots, A_n)$ of finite sets in $\Z^n$, we can consider tuples $(H_1,\ldots,H_n)$ of tropical hypersurfaces defined by tuples of tropical polynomials $(\varphi_1, \ldots, \varphi_n)$ such that $\supp(\varphi_i)=A_i, 1\leqslant i\leqslant n,$ as points in the space 
$$\EuScript{M}=\prod_1^n \T^{\vert A_i\vert}.$$

By $\EuScript{S}_0\subset\EuScript{M}$ denote the set of all transverse tuples of hypersurfaces. Obviously, the set $\EuScript{S}_0$ is open and everywhere dense. 
\begin{utver}\label{constS0}
The 2-intersection number $\iota_2(H_1,\ldots, H_n;\zeta)$ defined in $\ref{2-index}$ is constant on the connected components of $\EuScript{S}_0$. 
\end{utver}
$\vartriangleleft$
Suppose that the points $(H_1,\ldots, H_n)$ and $(H'_1,\ldots, H'_n)$ belong to the same connected component of $\EuScript{S}_0$. Then, there exists a one-to-one correspondence between the sets ${H_1\cap H_2\cap \ldots \cap H_n}$ and $H'_1\cap H'_2\cap \ldots \cap H'_n$, which maps every  $x\in H_1\cap H_2\cap \ldots \cap H_n$ to the point $x'\in H'_1\cap H'_2\cap \ldots \cap H'_n$ such that $\newton_x(H_i)=\newton_{x'}(H'_i)$ for every $i\in\{1,\ldots, n\}$. Therefore, the corresponding $2$-determinants in the right-hand side of (\ref{eq2ind}) coincide, which implies the sought equality 
$\iota_2(H_1,\ldots, H_n;\zeta)=\iota_2(H'_1,\ldots, H'_n;\zeta)$. 
$\vartriangleright$

The next step is to construct a codimension 1 set $\EuScript{S}_1\subset\EuScript{M}\setminus\EuScript{S}_0$ of ``almost transverse'' tuples of tropical hypersurfaces such that passing from one connected component of $\EuScript{S}_0$ to another through the points of $\EuScript{S}_1$ does not change the $2$-intersection number, and $\codim(\EuScript{M}\setminus(\EuScript{S}_0\cup\EuScript{S}_1))\geqslant 2$. Theorem \ref{main} will then follow, since the complement to a codimension 2 subset is connected. 

We shall prove below that the pieces of the sought set $\EuScript{S}_1$ are in one--to-- one correspondence with certain combinatorial structures of the form $(I_i,\, i\in I)$, $I\subset\{1,\ldots,n\}, I_i\subset A_i$, that we call elementary obstacles. 

\begin{defin} 
A collection of numbers $I\subset\{1,\ldots,n\}$ and pairs of points $I_i\subset A_i,\, i\in I$, is called an {\it elementary obstacle of type 1}, if the convex hull of the Minkowski sum $\sum_{i\in I} I_i$ has dimension $|I|-1$, and no proper subcollection $I'\subset I,\, I_i,\, i\in I'$, is an elementary obstacle of type 1.

A collection of numbers $I\subset\{1,\ldots,n\}$, a triple of points $I_j\in A_j,\, j\in I$, and pairs of points $I_i\subset A_i,\, i\in I\setminus\{j\}$, is called an {\it elementary obstacle of type 2}, if the convex hull of $I_j$ is a triangle, the convex hull of the Minkowski sum $\sum_{i\in I} I_i$ has dimension $|I|$, and no proper subcollection $I'\subset I,\, I_i,\, i\in I'$, is an elementary obstacle of type 1 or 2.

A one element collection $I=\{i\}$ and a triple of points $I_i\subset A_i$ is called an {\it elementary obstacle of type 3}, if the convex hull of $I_i$ is a segment.
\end{defin}

Below are shown all possible elementary obstacles in dimension 2.

\begin{center}
\begin{tikzpicture}[scale=0.8]
\draw[thick, dashed] (-2,-1)--(-2,1);
\draw[thick, blue, dashed] (-0.5,-1)--(-0.5,1);
\draw[fill] (-2,1) circle[radius=0.15];
\draw[fill] (-2,-1) circle[radius=0.15];
\draw[fill, blue] (-0.5,1) circle[radius=0.15];
\draw[fill, blue] (-0.5,-1) circle[radius=0.15];
\draw[thick, dashed] (3,-1)--(3,1)--(5,-1)--(3,-1);
\draw[fill] (3,1) circle[radius=0.15];
\draw[fill] (3,-1) circle[radius=0.15];
\draw[fill] (5,-1) circle[radius=0.15];
\draw[thick, dashed, blue] (6.5,-1)--(8.5,1);
\draw[fill, blue] (6.5,-1) circle[radius=0.15];
\draw[fill, blue] (8.5,1) circle[radius=0.15];
\draw[thick, dashed, blue] (12,0)--(16,0);
\draw[fill, blue] (12,0) circle[radius=0.15];
\draw[fill, blue] (14,0) circle[radius=0.15];
\draw[fill, blue] (16,0) circle[radius=0.15];
\end{tikzpicture}

\vspace{2ex}

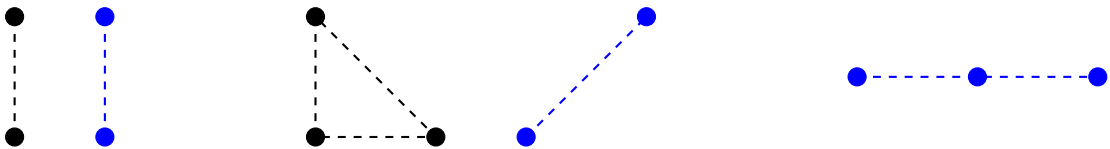
Figure 1. Elementary obstacles of types 1, 2 and 3 in dimension 2.
\end{center}

We say that tropical hypersurfaces $H_1,\ldots,H_n$ defined by tropical polynomials $\varphi_1,\ldots,\varphi_n$ have a non-transversality of type $k$, if there exists an elementary obstacle $I_i\subset A_i,\, i\in I$, of type $k$, such that for some point $x\in H_1\cap\ldots\cap H_n$ we have $\supp_x(\varphi_i)=I_i,\, i\in I$.

The rest of Section 2 will be spent to observe that generic tuples of hypersurfaces with a non-transversality of one of the three types form the sought set $\EuScript{S}_1$. In other words, if we travel between two tuples of transversal hypersurfaces along a generic path in $\EuScript{M}$, then we shall encounter finitely many generic tuples with a non-transversality of type $k$, and passing through them will not change the 2-intersection number. Thus, passing through generic tuples with a non-transversality of type $k$ plays the role of Reidemeister moves in knot theory. The figure below shows all such ``Reidemeister moves'' in dimension two.

\begin{center}
\begin{tikzpicture}[scale=0.9]
\draw[ultra thick] (-2,0)--(0,0)--(0,-2);
\draw[ultra thick] (0,0)--(2,2);\node[below right] at (1.5,1.5) { $C_1$};
\draw[ultra thick][blue](-1,2)--(0,1)--(2,-1);\node[above right, blue] at (-1.5,1.2) { $C_2$};
\draw[fill, red] (0.5,0.5) circle[radius=0.1];
\draw [ultra thick, blue] (-3,2)--(-4,0)--(-3,-2);
\draw [ultra thick, blue] (-4,0)--(-8.5,0);
\draw [ultra thick] (-4.5,2.5)--(-5.5,1)--(-4.8,-0.8);
\draw [ultra thick] (-5.5,1)--(-7,1);
\draw [ultra thick] (-8,2.5)--(-7,1)--(-8,-1);
\draw [ultra thick, blue] (-9.5,1.5)--(-8.5,0)--(-9.5,-2);
\draw [fill, red] (-5.13,0) circle[radius=0.1];
\draw [fill, red] (-7.5,0) circle[radius=0.1];
\node[right, blue] at (-4,1.7) {$C_2$};
\node[left] at (-7.7, 2.2) {$C_1$};
\draw[ultra thick] (4,0)--(8,0);
\draw[line width=0.12 cm][blue] (6,-2)--(6,2);
\draw[fill, red] (6,0) circle[radius=0.1];
\node[right, blue] at (6,1.5) {$C_2$};
\node[above] at (4.5,0) {$C_1$};
\end{tikzpicture}
\end{center}

\begin{center}
\begin{tikzpicture}[scale=0.9]
\draw[ultra thick] (-2,0)--(0,0)--(0,-2);
\draw[ultra thick] (0,0)--(2,2);\node[below right] at (1.5,1.5) { $C_1$};
\draw[ultra thick][blue](-2,2)--(0,0)--(2,-2);\node[above right, blue] at (-1.5,1.5) { $C_2$};
\draw[fill, red] (0,0) circle[radius=0.1];
\draw [ultra thick, blue] (-3,2)--(-4,0)--(-3,-2);
\draw [ultra thick, blue] (-4,0)--(-5.5,0);
\draw [ultra thick] (-4.5,1.5)--(-5.5,0)--(-5,-1.7);
\draw [ultra thick, red] (-5.5,0)--(-7,0);
\draw [ultra thick] (-8,1.5)--(-7,0)--(-7.7,-1.7);
\draw [ultra thick, blue] (-7,0)--(-8.5,0);
\draw [ultra thick, blue] (-9.5,1.5)--(-8.5,0)--(-9.5,-2);
\draw [fill, red] (-5.5,0) circle[radius=0.1];
\draw [fill, red] (-7,0) circle[radius=0.1];
\node[right, blue] at (-4,1.7) {$C_2$};
\node[left] at (-8, 1.5) {$C_1$};
\draw[ultra thick] (4,0)--(8,0);
\draw[line width=0.12 cm][blue] (6,-2)--(6,2);
\draw[fill, red] (6,0) circle[radius=0.1];
\node[right, blue] at (6,1.5) {$C_2$};
\node[above] at (4.5,0) {$C_1$};
\end{tikzpicture}
\end{center}

\begin{center}
\begin{tikzpicture}[scale=0.9]
\draw[ultra thick] (-2,0)--(0,0)--(0,-2);
\draw[ultra thick] (0,0)--(2,2);\node[below right] at (1.5,1.5) { $C_1$};
\draw[ultra thick][blue](-2,1)--(0,-1)--(1,-2);\node[above right, blue] at (-2,1) { $C_2$};
\draw[fill, red] (0,-1) circle[radius=0.1];
\draw[fill, red] (-1,0) circle[radius=0.1];
\draw [ultra thick, blue] (-3,2)--(-4,0)--(-3,-2);
\draw [ultra thick, blue] (-4,0)--(-8.5,0);
\draw [ultra thick] (-4.5,0.8)--(-5.5,-0.7)--(-5,-2.4);
\draw [ultra thick] (-5.5,-0.7)--(-7,-0.7);
\draw [ultra thick] (-8,0.8)--(-7,-0.7)--(-7.7,-2.4);
\draw [ultra thick, blue] (-9.5,1.5)--(-8.5,0)--(-9.5,-2);
\draw [fill, red] (-5.075,0) circle[radius=0.1];
\draw [fill, red] (-7.43,0) circle[radius=0.1];
\node[right, blue] at (-4,1.7) {$C_2$};
\node[left] at (-7.9, 0.5) {$C_1$};
\draw[ultra thick] (4,0)--(8,0);
\draw[ultra thick][blue] (5,-2)--(5,2);
\draw[ultra thick][blue] (7,-2)--(7,2);
\draw[fill, red] (5,0) circle[radius=0.1];
\draw[fill, red] (7,0) circle[radius=0.1];
\node[right, blue] at (7,1.5) {$C_2$};
\node[above] at (4.5,0) {$C_1$};
\end{tikzpicture}

\vspace{2ex}

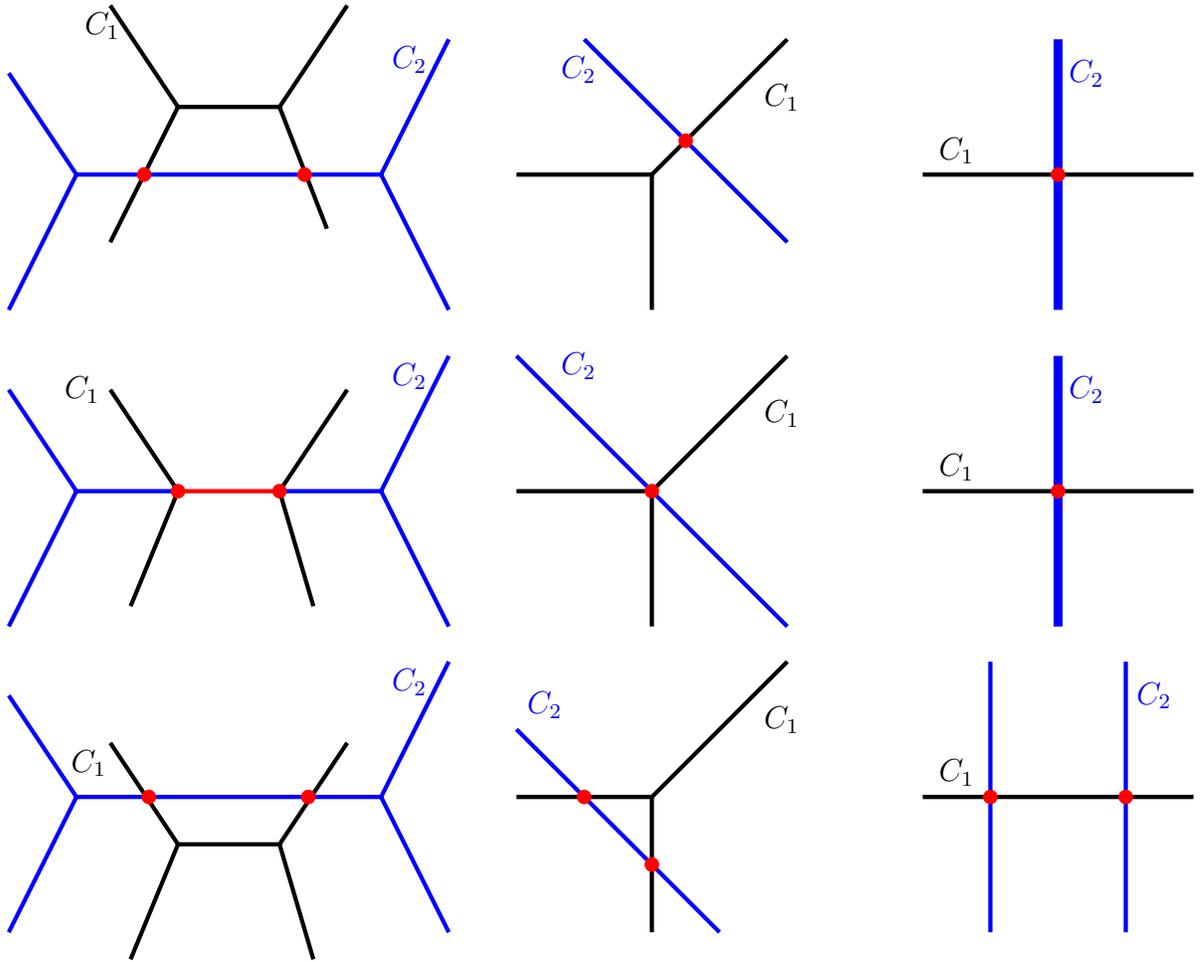
Figure 2. ``Reidemeister moves'' of types 1, 2 and 3 in dimension 2.
\end{center}

\subsection{Building the Walls}\label{buildwall}

\begin{rem}
Here we consider tuples of finite sets in $\Z^n$. Thus the words ``simplex'' and ``interval'' mean a set of all vertices of a simplex or an interval, respectively. 
\end{rem}

In the notation of Theorem \ref{main}, consider a $n$-tuple $A=(A_1,\ldots, A_n)$ of finite sets in $\Z^n$ such that $\conv(A_i)=P_i, 1\leqslant i\leqslant n$ and the polytopes $P_1,\ldots, P_n$ are 2-developed with respect to a point $\zeta\in\Z^n$. Let $B=(B_1,\ldots, B_n)$ be an arbitrary subtuple $B=(B_1,\ldots, B_n), B_i\subset A_i$ of simplices. 

\begin{defin}\cite{sturmfels}\label{codim}
For an arbitrary subset $I\subset\{0,\ldots,n\}$, we define its {\it codimension}, which we denote by $\codim(I)$, as follows:
$$ \codim(I)=\dim(\conv(\sum_{i\in I} B_i))-|I|. $$
The {\it codimension of the tuple} $B$ is defined by the following equality: 
$$\codim(B)=\min_{I\subset\{0,\ldots,n\}}\codim(I).$$
\end{defin}

Let $\varphi=(\varphi_1,\ldots, \varphi_n)$ be a $n$-tuple of tropical polynomials, where ${\varphi_i=\bigoplus_{a\in A_i} c_{a,i}\odot x^{\odot a}}$. 
To every point $b\in B_i$ one can associate the point $b^{\varphi}=(b,c_{b,i})\in\Z^n \times\T$. We denote by $B^{\varphi}$ the tuple $(B_1^{\varphi_1}, \ldots, B_n^{\varphi_n})$, where $B_i^{\varphi_i}=\{b^{\varphi_i}\mid b\in B_i\}$.

\begin{defin}
In the previous notation, given a tuple $B$, by ${L_B}$ we denote the affine subspace consisting of all the points $\varphi\in\EuScript{M}$ such that the following equality holds:
$\dim(\conv(\sum_1^n B_i^{\varphi_i}))=\dim(\conv(\sum_1^n B_i))$. 

A subtuple $B\subset A$ is called a {\it trouble}, if $L_B\neq \EuScript{M}$ and $|B_i|>1$ for all $i$. We say that a  trouble $B$ is an {\it obstacle}, if $\codim(L_B)=1$. In this case, we call the hyperplane  $L_B$ a {\it wall}. 
\end{defin}

The sets $L_B$ corresponding to each of the troubles $B\subset A$ cover the set $\EuScript{M}\setminus\EuScript{S}_0$. In these terms, the sought codimension 1 set $\EuScript{S}_1\subset\EuScript{M}$ is a set, which contains all the points $\varphi=(\varphi_1,\ldots, \varphi_n)$ such that all the walls containing $\varphi$ coincide and for every trouble $B\subset A$, $\codim(L_B)>1$ implies $\varphi\notin L_B$. In order to obtain the explicit description of the set $\EuScript{S}_1$, we need first to describe and classify the obstacles $B\subset A$. 

Let $B\subset A$ be a trouble. It is obvious, that $\codim(B)\leqslant 0$. Consider the following cases: 
\begin{enumerate}
\item $\codim(B)\leqslant -1$;
\item $\codim(B)=0$.
\end{enumerate}

{\bf Case 1.} Consider a subtuple $B'\subset B$ of intervals $B'_i\subset B_i$. For every set ${I\subset\{0,\ldots,n\}}$, the inequality $\dim(\conv(\sum_{i\in I} B'_i))\leqslant \dim(\conv(\sum_{i\in I} B_i))$ holds, therefore $\codim(B')\leqslant\codim(B)$. Let $I_{\min}$ be the minimal subset $I\subset \{1,\ldots, n\}$ such that $\codim(B')=\codim(I)$. 

\begin{rem}
The set $I_{\min}$ is non-empty, because $\codim(\varnothing)=0$,\\ while $\codim(I_{\min})\leqslant -1$.
\end{rem}

By $B''$ we denote the subtuple $(B'_i, i\in I_{\min})$. 
\begin{rem}
It is easy to show that $\codim(B')=\codim(B'')$. We obviously have $\codim(B'')\leqslant\codim(B')$. Assume that $\codim(B'')<\codim(B')$. Then, there exists a subset $J\subsetneq I_{\min}$ of such that $\codim(J)<\codim(I_{\min})$, which contradicts with the choice of the subset $I_{\min}$.
\end{rem}
\begin{utver}\label{codim1}
The following equality holds: $\codim(L_{B''})=-\codim(B'')$.
\end{utver}

$\vartriangleleft$ Without loss of generality, suppose that $I_{min}={1,\ldots,m}$ for some $m\leqslant n$. Then the tuple $B''$ is a $m$-tuple $B'_1,\ldots, B'_m$ of intervals in $\Z^n$. 
For every $1\leqslant i\leqslant m$, by $v_i$ we denote the vector $\overrightarrow{B'_i}\in\Z^n$. By $v_{ij}$ we denote the $j$-th component of the coordinate vector of $v_i$ in the standard basis. Consider an arbitrary point $\varphi=(\varphi_1,\ldots, \varphi_n)\in\EuScript{M}$. Let $\alpha_i$ and $\beta_i$ be the coefficients of $\varphi_i$ corresponding to each of the endpoints of the interval $B'_i$. In these terms, the set $L_{B''}$ consists of the points $\varphi\in\EuScript{M}$ such that the following equality holds:
\begin{equation}\label{rankcodim1}
m+\codim(B'')=\rk \begin{pmatrix}
v_{11} & \dots & v_{1n} \\
\vdots & \ddots & \vdots \\
v_{m1} & \dots & v_{mn}
\end{pmatrix}=\rk \begin{pmatrix}
v_{11} & \dots & v_{1n} & (\beta_1-\alpha_1) \\
\vdots & \ddots & \vdots \\
v_{m1} & \dots & v_{mn} & (\beta_m-\alpha_m) 
\end{pmatrix}
\end{equation}
We can suppose without loss of generality that the first $m+\codim(B'')$ columns span the column space of the first matrix. Then the equality (\ref{rankcodim1}) means that the last column of the second matrix can be expressed as their linear combination. Thus the sought codimension of the plane $L_{B''}$ equals ${m-m-\codim(B'')=-\codim(B'')}$.
$\vartriangleright$

\begin{sledst}
If $\codim(B)\leqslant -2$, then $\codim(L_B)\geqslant 2$. 
\end{sledst}

\begin{proof}
Obviously, $L_B\subset L_{B''}$. Therefore,  $\codim(L_B)\geqslant \codim(L_{B''})$. Applying \ref{codim1}, we obtain ${\codim(L_B)\geqslant \codim(L_{B''})=-\codim(B'')\geqslant -\codim(B)\geqslant 2}$.
\end{proof}

\begin{sledst}
In the previous notation, if $B$ is an obstacle, and $\codim(B)=-1$, then $L_B=L_{B''}$. In this case, we call $B''$ an {\it elementary obstacle} corresponding to the obstacle $B$.
\end{sledst}

\begin{proof}
The equalities $\codim(B)=-1$ and $\codim(L_B)=1$ imply \\that $1\leqslant -\codim(B'')=\codim(L_{B''})\leqslant 1$. So, $\codim(L_{B''})=-\codim(B'')=1$. Therefore, $1=\codim(L_B)\geqslant\codim(L_{B''})=1$, which finishes the proof. 
\end{proof}

{\bf Case 2a.} Suppose $\codim(B)=0$ and the convex hull of at least one of $B_i$'s is not a segment. Then we can choose $B'\subset B$ to be a subtuple consisting of a triangle and $(n-1)$ intervals $B'_i\subset B_i$. Let $M_{B'}=\min\codim(I)$, where the minimum is taken over all the subsets $I\subset\{1,\ldots, n\}$ such that the tuple $(B'_i, i\in I)$ contains the triangle. By $I_{\min}$ denote the minimal set such that $\codim(I_{\min})=M_{B'}$ and the tuple $B''=(B'_i, i\in I_{\min})$ contains the triangle. 

\begin{utver}\label{codim2}
In the previous notation, the following equality holds: \\$\codim(L_{B''})=-M_{B'}+1.$
\end{utver}

$\vartriangleleft$
The proof is almost the same as the one of Proposition \ref{codim1}. The only difference is that the triangle gives rise to two vectors instead of one, thus, we will deal with $(m+1)\times n$ and $(m+1)\times (n+1)$-matrices of rank $m+M_{B'}$. Therefore, in this case, the sought codimension equals $m+1-m-M_{B'}=-M_{B'}+1$.
$\vartriangleright$

\begin{sledst}
If $\codim(B'')\leqslant -1$, then $\codim(L_B)\geqslant 2$.
\end{sledst}

$\vartriangleleft$
Using \ref{codim2}, we have \\ ${\codim(L_B)\geqslant\codim(L_{B''})
=-\codim(B'')+1\geqslant 2}$, which finishes the proof.
$\vartriangleright$

\begin{sledst}
If $B$ is an obstacle, and $\codim(B)=0$, $L_B=L_{B''}$. In this case, $B''$ is called an {\it elementary obstacle} corresponding to the obstacle $B$.
\end{sledst}

$\vartriangleleft$
The equalities $\codim(B)=0$ and $\codim(L_B)=1$ imply that\\ $1=-\codim(B)+1\leqslant-M_{B'}+1=\codim(L_{B''})\leqslant
\codim(L_B)=1$, which finishes the proof.

Moreover, from the proof, it follows that $M_{B'}=0$.
$\vartriangleright$

{\bf Case 2b.} Suppose $\codim(B)=0$ and the convex hull of each of $B_i$'s is not a segment. If $|B_i|=2$ for all $i$, then $\codim L_B=0$, and if $|B_i|>3$ for some $i$, then $\codim L_B>1$. So, if $B$ is an obstacle, then $|B_i|=3$ for some $i$. In this case we denote the one element subtuple $(B_i)$ by $B''$ and observe that $L_B=L_{B''}$ provided that $\codim L_B=1$.

\begin{defin}\label{elemobst}
An {\it elementary obstacle} is a subtuple $\mathsf K=(K_i\mid i\in I)$ of subsets $K_i\subset A_i$, where $I\subset\{1,\ldots,n\}$, belonging to one of the following types:
\begin{description}
\item{\bf Type 1:} A codimension $-1$ tuple of intervals such that $I_{\min}=I$;
\item{\bf Type 2:} A codimension $0$ tuple consisting of one triangle and $|I|-1$ intervals which satisfies the following properties:
\begin{itemize}
\item $I_{\min}=I$;
\item there exists no elementary obstacle $\mathsf M=(M_j\mid j\in J\subset I)$ of Type 1 such that $ M_j\subset K_j$ for every $j\in J$.
\end{itemize} 
\item{\bf Type 3:} $I=\{i\}$, and $K_i$ consists of three points on a line.
\end{description}
\end{defin}

Let $\mathsf{K}$ be an elementary obstacle. Without loss of generality, suppose that $I=\{1,\ldots, m\}$ for some $m\leqslant n$. If $\mathsf{K}$ is of the first type, then for every $1\leqslant i\leqslant m$, by $v_i$ we denote the vector $\overrightarrow{K_i}\in\Z^n$. Otherwise, suppose that $K_1$ is a triangle. Then we set $v_0=\overrightarrow{K'_0}$, and $v_1=\overrightarrow{K'_1}$ for any two edges $K'_0$ and $K'_1$ of $K_1$ and $v_i=\overrightarrow{K_i}$ for every $2\leqslant i\leqslant m$.

By $v_{ij}$ we denote the $j$-th component of the coordinate vector of $v_i$ in the standard basis. Consider an arbitrary point $\varphi=(\varphi_1,\ldots, \varphi_n)\in\EuScript{M}$. Let $\alpha_i$ and $\beta_i$ be the coefficients of $\varphi_i$ corresponding to each of the endpoints of the intervals $K_i$ and $K'_i$. 

In these terms, if the obstacle $\mathsf{K}$ is of the first type, then the wall  $L_{\mathsf{K}}$ consists of the points $\varphi\in\EuScript{M}$ such that the following equality holds:
\begin{equation*}
m-1=\rk \begin{pmatrix}
v_{11} & \dots & v_{1n} \\
\vdots & \ddots & \vdots \\
v_{m1} & \dots & v_{mn}
\end{pmatrix}=\rk \begin{pmatrix}
v_{11} & \dots & v_{1n} & (\beta_1-\alpha_1) \\
\vdots & \ddots & \vdots \\
v_{m1} & \dots & v_{mn} & (\beta_m-\alpha_m) 
\end{pmatrix}
\end{equation*}

Therefore, $L_{\mathsf{K}}$ is defined by the equations of the following form:
\begin{equation}\label{defeq1}
\det \begin{pmatrix}
v_{1{j_1}} & \dots & v_{1{j_{m-1}}} & (\beta_1-\alpha_1) \\
\vdots & \ddots & \vdots \\
v_{m{j_1}} & \dots & v_{m{j_{m-1}}} & (\beta_m-\alpha_m) 
\end{pmatrix}=0
\end{equation}

Thus the defining equations of $L_{\mathsf{K}}$ are linear and employ only the coefficients $\alpha_i$ and $\beta_i$ as variables. Moreover, for every $1\leqslant i\leqslant m$, $\alpha_i$ and $\beta_i$ occur at least in one of the defining equations (\ref{defeq1}) with nonzero coefficients, since otherwise, we would find a $(m-1)$-tuple of linearly dependent vectors $(v_1,\ldots,\hat{v_i}, v_m)$, which would mean that $\mathsf{K}$ is not an elementary obstacle. 

The same arguments work for the case of $\mathsf{K}$ being an elementary obstacle of type 2. In this case, the wall $L_{\mathsf{K}}$ consists of the points $\varphi\in\EuScript{M}$ such that the following equality holds:

\begin{equation*}
m=\rk \begin{pmatrix}
v_{01} & \dots & v_{0n} \\
\vdots & \ddots & \vdots \\
v_{m1} & \dots & v_{mn}
\end{pmatrix}=\rk \begin{pmatrix}
v_{01} & \dots & v_{0n} & (\beta_0-\alpha_0) \\
\vdots & \ddots & \vdots \\
v_{m1} & \dots & v_{mn} & (\beta_m-\alpha_m) 
\end{pmatrix}
\end{equation*}

Therefore, the wall $L_{\mathsf{K}}$ is defined by the equations of the following form:
\begin{equation}\label{defeq2}
\det \begin{pmatrix}
v_{0{j_1}} & \dots & v_{0{j_{m}}} & (\beta_0-\alpha_0) \\
\vdots & \ddots & \vdots \\
v_{0{j_1}} & \dots & v_{m{j_{m}}} & (\beta_m-\alpha_m) 
\end{pmatrix}=0
\end{equation}
The same arguments as above imply that the equations \ref{defeq2} are linear and employ only the coefficients $\alpha_i$ and $\beta_i$ as variables, moreover, each of the coefficients $\alpha_i$ and $\beta_i$ occurs at least in one of the equations with a nonzero coefficient. 

The same is obviously valid for elementary obstacles of type 3, so, we have proved the following

\begin{lemma}\label{elementobst}
In the previous notation, let $\mathsf K_1$ and $\mathsf K_2$ be elementary obstacles. Then $L_{\mathsf K_1}=L_{\mathsf K_2}$ if and only if $\mathsf K_1=\mathsf K_2$.
\end{lemma}

The following statements easily follow from Lemma \ref{elementobst}. 

\begin{sledst}
Each obstacle $B\subset A$ has a unique elementary obstacle~${K\subset B}$.
\end{sledst}
\begin{sledst}
Let $B_1,B_2\subset A$ be obstacles with the elementary obstacles $K_1$ and $K_2$ respectively. Then, the walls $L_{B_1}$ and $L_{B_2}$ coincide if and only if $K_1=K_2$.
\end{sledst}
 
\begin{defin}
The set $\EuScript{S}_1$ of {\it almost transverse} tropical hyperplanes is defined to be the set of all $\varphi=(\varphi_1,\ldots,\varphi_n)\in \EuScript{M}\setminus\EuScript{S}_0$ such that the following conditions are satisfied:
\begin{itemize}
\item there exists a wall containing $\varphi$;
\item all the walls containing $\varphi$ coincide;
\item for every trouble $B$, $\codim(L_B)\geqslant 2$ implies that $\varphi\not\in L_B$.
\end{itemize} 
\end{defin}

\begin{sledst}
For every point $\varphi\in\EuScript{S}_1$ there exists a unique elementary obstacle $\mathsf{K}$ such that $\varphi\in L_{\mathsf{K}}$. Moreover, in a small neighborhood of the point $\varphi$, we have $\EuScript{S}_1=L_{\mathsf K}$.
\end{sledst}

\begin{sledst}
The following inequality holds: $\codim(\EuScript{M}\setminus(\EuScript{S}_0\cup\EuScript{S}_1))\geqslant 2$.
\end{sledst}

\begin{sledst}
The set $\EuScript{S}_0\cup\EuScript{S}_1$ is connected.
\end{sledst}

We now explicitly describe every tuple $\varphi=(\varphi_1,\ldots,\varphi_n)\in\EuScript{S}_1$.

\begin{defin}
A {\it generic extension} of an elementary obstacle ${\mathsf K}=(K_i,\, i\in I)$ is a collection ${\mathsf K}'=(K'_1,\ldots,K'_n),\, K'_i\subset A_i$, satisfying the following properties:
\begin{itemize}
\item Assume that ${\mathsf K}$ is of type 1. Then there are three possibilities for $K'$
\subitem (1) $K'_i=K_i$ for $i\in I$, and otherwise $K'_i$ is a pair of points such that the convex hull of the Minkowski sum $\sum_i K'_i$ has codimension 1.
\subitem (2) The same as (1), but $K'_i\supset K_i$ is a triangle for one $i\in I$, and the Minkowski sum $\sum_i K'_i$ is not contained in a hyperplane. 
\subitem (3) The same as (1), but $K'_i$ is a triangle for one $i\notin I$, and the Minkowski sum $\sum_i K'_i$ is not contained in a hyperplane. 
\item Assume that ${\mathsf K}$ is of type 2. Then $K'_i=K_i$ for $i\in I$, and otherwise $K'_i$ is a pair of points such that  the Minkowski sum $\sum_i K'_i$ is not contained in an affine hyperplane.
\item Assume that ${\mathsf K}$ is of type 3. Then $K'_i=K_i$ for $i\in I$, and otherwise $K'_i$ is a pair of points such that  the Minkowski sum $\sum_i K'_i$ is not contained in an affine hyperplane.
\end{itemize}
\end{defin}

\begin{sledst}\label{badinters}
According to the classification of the elementary obstacles $\mathsf{K}\subset A$ given in \ref{elemobst}, each point $\varphi=(\varphi_1,\ldots,\varphi_n)\in\EuScript{S}_1\cap L_{\mathsf K}$ defines a tuple $H=(H_1,\ldots, H_n)$ of tropical hypersurfaces that intersects as follows: at every non-transversal point $x$ of the intersection $H_1\cap\ldots\cap H_n$, the local support sets $\supp_x(\varphi_i)$ form a generic extension of ${\mathsf K}$.
\end{sledst}

In particular, if the elementary obstacle ${\mathsf K}=(K_i,\, i\in I)$ is of 
\begin{description}
\item{\bf type 1:} the intersection $C$ of the $|I|$ facets $C_i$ of the tropical hypersurfaces $(H_i, i\in I)$ that are dual to the intervals $(K_i, i\in I)$ is ${(n-|I|+1)}$-dimensional. The intersection of $C$ with the rest of $H_i,~i\notin I$, is a graph with vertices of degree 3 and 1. All the other intersections of the  hypersurfaces $(H_1,\ldots,H_n)$ are transverse; 
\item{\bf type 2:} {\bf exactly} one of the intersection points $x\in H_1\cap\ldots\cap H_n$ is not smooth, since $\newton_x(H_i)$ is a triangle for some $i\in I$. All the other intersections of the hypersurfaces $(H_1,\ldots,H_n)$ are transverse;
\item{\bf type 3:} The intersection $H_1\cap\ldots\cap H_n$ is transversal in the sense that at each of its finitely many points $x$ the local Newton polytopes $\newton_x(H_i)$ are transversal segments, but at some of these points the support set $\supp_x(\varphi_i),\, i\in I,$ consists of three points instead of two. 
\end{description}

Our next step is to show that the 2-intersection number $\iota_2$ is constant on  $\EuScript{S}_0\cup\EuScript{S}_1$. 

\subsection{Passing through the walls}\label{passwall}
\begin{defin}
We say that a point $x\in\T^n\times\{(H_1,\ldots,H_n)\}\subset\T^n\times\EuScript{M}$~belongs to the {\it stable intersection} of tropical hypersurfaces $H_1,\ldots,H_n$ (denote $H_1\cap_{st}\ldots\cap_{st}H_n$), if for every $\varepsilon>0$ there exists $\delta>0$ such that for any tuple of translations $(T_{\overrightarrow{v_1}},\ldots, T_{\overrightarrow{v_n}})$, where $|{\overrightarrow{v_i}}|<\delta$, there exists a point $x'\in T_{\overrightarrow{v_1}}(H_1)\cap\ldots\cap T_{\overrightarrow{v_n}}(H_n)$ such that $|x-x'|<\varepsilon$.
\end{defin}
\begin{rem}
If $H_1\pitchfork\ldots\pitchfork H_n$, then $H_1\cap_{st}\ldots\cap_{st}H_n =H_1\cap\ldots\cap H_n$. 
\end{rem}
\begin{lemma}\label{epsilon}
For every open $U\subset\T^n$ such that ${U\times\{(H_1,\ldots,H_n)\}\subset\T^n\times\EuScript{M}}$ contains the stable intersection $H_1\cap_{st}\ldots\cap_{st}H_n$, there exists an open ${V\subset\EuScript{M}}$ satisfying the following properties: $(H_1,\ldots,H_n)\in V$ and for every ${(H'_1,\ldots,H'_n)\in V}$ the stable intersection ${H'_1\cap_{st}\ldots\cap_{st}H'_n}$ is contained in the set ${U\times\{(H'_1,\ldots,H'_n)\}}$.
\end{lemma}
$\vartriangleleft$ Consider the projection $\pi\colon\T^n\times\EuScript{M}\to\EuScript{M}$. It suffices to show that for any element $\varphi\in\EuScript{M}$ and for any open $\pi^{-1}(\varphi)\subset U\subset\T^n\times\{\varphi\}$ there exists an open $V\subset\EuScript{M}$ such that $\varphi\in V$ and for every $\varphi'\in V$ the preimage $\pi^{-1}(\varphi')$ is contained in $U\times V$. Choose an arbitrary neighbourhood $I\ni\varphi$ and consider the set $K= \pi(\overline{U\times I}\setminus\pi^{-1}(I))$. This set is compact and does not contain the point $\varphi$, so there exists a neighbourhood $V\ni \varphi$ such that $K\cap V=\varnothing$, which finishes the proof of the lemma.
$\vartriangleright$

Consider a point $(H_1,\ldots,H_n)\in \EuScript{S}_1$. Besides the non-transverse intersection points mentioned in \ref{badinters}, the hypersurfaces $H_1,\ldots,H_n$ have finitely many transverse intersection points $Q_1, \ldots, Q_l$. The following statement is obvious.
\begin{utver}\label{goodpoints}
For every transverse intersection $Q_j\in H_1\cap\ldots \cap H_n$, there exist an open ${U_j\ni Q_j}$ and an open ${V_j\ni (H_1,\ldots,H_n)}$ such that for every ${(H'_1,\ldots, H'_n)\in V_j}$ all the intersection points of the hypersurfaces $H'_1,\ldots,H'_n$ are transverse in $U_j$.  
\end{utver} 

\begin{theor}\label{calculS1}
For every $(H_1,\ldots,H_n)\in\EuScript{S}_1$, $\iota_2$ is constant on ${V_0=\bigcap_{j=1}^m (V_j\cap V)\subset\EuScript{M}}$.
\end{theor}

\begin{proof}
In Subsection \ref{buildwall}, we obtained the following classification of the walls in $\EuScript M$. Every wall is the hyperplane $L_{\mathsf K}$ corresponding to an elementary obstacle $\mathsf K\subset A$, i.e., a subtuple $\mathsf K=(K_i\mid i\in I)$, where $K_i\subset A_i$ and $I\subset\{1,\ldots,n\}$, belonging to one of the following types:
\begin{description}
\item{\bf Type 1:} A codimension $-1$ tuple of intervals such that $I_{\min}=I$;
\item{\bf Type 2:} A codimension $0$ tuple consisting of one triangle and $|I|-1$ intervals such that $I_{\min}=I$.
\item{\bf Type 3:} One triple of points on a line.
\end{description}

For each of the types of walls $L_{\mathsf K}$, we will show that passing through points $\varphi\in L_{\mathsf K}$ does not change the 2-intersection number. 
Namely, if connected components $S$ and $S'$ of the set $\EuScript{S}_0$ are  separated by a wall $L_{\mathsf K}$, then for any points $(H_1,\ldots, H_n)\in S$ and $(H'_1,\ldots, H'_n)\in S'$, there exists a path $\gamma\colon [0,1]\to \EuScript{S}_0\cup\EuScript{S}_1$ with the endpoints $(H_1,\ldots, H_n)$ and $(H'_1,\ldots, H'_n)$ such that $\iota_2(\gamma(0);\zeta)=\iota_2(\gamma(1);\zeta)$.

{\bf Case 1.} Consider a point $(H_1,\ldots, H_n)\in\EuScript{S}_1$ that belongs to a wall $L_{\mathsf K}$ corresponding to the elementary obstacle $\mathsf K=(K_i\mid i\in I)$ of type 1. Moreover, we may assume that the point $(H_1,\ldots, H_n)$ is generic in $L_{\mathsf K}$.

Without loss of generality, assume that $I=\{1,\ldots,m\}$ for some $m\leqslant n$. By $C_i$ we denote the facet of $H_i$ which is dual to the interval $K_i$, and by $C$ we denote the intersection $C_1\cap\ldots\cap C_m$.

In these terms, for every point $x\in C\cap H_{m+1}\cap\ldots\cap H_n$, the tuple $\mathsf K$ is the elementary obstacle corresponding to the obstacle $B_x=(\newton_x(H_1),\ldots, \newton_x(H_n))$.
It immediately follows from Corollary \ref{badinters} that for every $i\in I$, we have \\ ${H_1\pitchfork\ldots\pitchfork H_{i-1}\pitchfork\hat{H_i}\pitchfork H_{i+1}\pitchfork\ldots\pitchfork H_n}$ (where $\hat{H_i}$ means that the hypersurface $H_i$ is omitted), i.e., these hypersurfaces intersect in a tropical curve. Without loss of generality, we take $i=1$, and by $\Sigma$ we denote the tropical curve $\Sigma=H_2\cap\ldots\cap H_n$. 

In the previous notation, we have $\Sigma\cap H_1=(\Sigma\cap C)\cup(Q_1\cup\ldots\cup Q_l)$, where $Q_j$ are the transverse intersection points of the hypersurfaces $H_1,\ldots, H_n$. Pick a normal $\overrightarrow{w}$ to the hyperplane containing the facet $C_1$ of the tropical hypersurface $H_1$. Thus we obtain a family of hypersurfaces $H_1(\varepsilon)=T_{\varepsilon\overrightarrow{w}}(H_1)$ (where $T_{\overrightarrow{v}}$ stands for the shift by a  vector $\overrightarrow{v}$) parametrized by $\varepsilon\in\R$. By $C_1(\varepsilon)$ we denote the image of the facet $C_1$ of the hypersurface $H_1$.

The neighbourhood $V_0\ni (H_1,\ldots, H_n)$ constructed above contains an open ball $\EuScript B\ni (H_1,\ldots, H_n)$, therefore, there exists $\varepsilon_0>0$ such that for every $\varepsilon \in (-\varepsilon_0, \varepsilon_0)$ the point $(H_1(\varepsilon), H_2,\ldots, H_n)$ belongs to $\EuScript B$. Fix an arbitrary ${0<\varepsilon_1<\varepsilon_0}$. For every $\varepsilon\neq 0$ in $[-\varepsilon_1, \varepsilon_1]$, we have $H_1(\varepsilon)\pitchfork H_2\pitchfork\ldots\pitchfork H_n$. So, our next step is to compare the 2-intersection numbers $\iota_2(H_1(-\varepsilon_1), H_2, \ldots, H_n;\zeta)$ and $\iota_2(H_1(\varepsilon_1), H_2, \ldots, H_n;\zeta)$ for the tuples $H_1(\pm\varepsilon_1), H_2, \ldots, H_n$ which are contained in different connected components of the set $\mathsf S_0$ separated by the wall $L_{\mathsf K}$. 

From Proposition \ref{goodpoints}, it follows that for every transverse intersection point \\${Q_j\in H_1\cap H_2\cap\ldots\cap H_n}$, there exists a neighbourhood ${U_j\ni Q_j}$, such that for every~${0\neq\varepsilon\in [-\varepsilon_1, \varepsilon_1]}$, the neighbourhood $U_j$ contains a unique transverse intersection point~${Q_j(\varepsilon)\in H_1(\varepsilon)\cap H_2\cap\ldots\cap H_n}$. Therefore, for every ${0\neq\varepsilon\in [-\varepsilon_1, \varepsilon_1]}$, we have: 
\begin{multline}\label{iotasum}
\iota_2(H_1(\varepsilon), \newton_{x}(H_2),\ldots, H_n;\zeta)= \sum\limits_{x\in C_1(\varepsilon)\cap\Sigma}{\mathrm{det}_2}(\newton_x(H_1(\varepsilon)), H_2,\ldots \newton_x(H_n),\zeta) \\+\sum_{j=1}^{j=l}{\mathrm{det}_2}(\newton_{Q_j(\varepsilon)}(H_1(\varepsilon)), \newton_{Q_j(\varepsilon)}(H_2),\ldots,  \newton_{Q_j(\varepsilon)}(H_n),\zeta). 
\end{multline}

Moreover, from the construction of the interval $[-\varepsilon_1, \varepsilon_1]$, it follows that in the summands ${\mathrm{det}_2}(\newton_{Q_j(-\varepsilon_1)}(H_1(-\varepsilon_1)),\ldots \newton_{Q_j(-\varepsilon_1)}(H_n),\zeta)$ \\ and ${\mathrm{det}_2}(\newton_{Q_j(\varepsilon_1)}(H_1(\varepsilon_1)),\ldots \newton_{Q_j(\varepsilon_1)}(H_n),\zeta)$ from the right-hand side of (\ref{iotasum}), the same Newton intervals occur, thus the two sums over the corresponding transverse intersection points coincide. Therefore, it suffices to deduce the following equality:
\begin{multline}
\sum\limits_{x\in C_1(-\varepsilon_1)\cap\Sigma}{\mathrm{det}_2}(\newton_x(H_1(-\varepsilon_1)), \newton_{x}(H_2),\ldots \newton_x(H_n),\zeta)=\\ \sum\limits_{x\in C_1(\varepsilon_1)\cap\Sigma}{\mathrm{det}_2}(\newton_x(H_1(\varepsilon_1)), \newton_{x}(H_2),\ldots \newton_x(H_n),\zeta)
\end{multline}

\begin{rem}\label{genericity}
Recall that the point $\varphi=(H_1,\ldots, H_n)\in L_{\mathsf K}$ is chosen to be generic (i.e., the corresponding dual subdivisions of the Newton polytopes $\newton(H_i)$ are simplicial for all $i\in\{1,\ldots,n\}$), thus, we can also assume that it satisfies the following conditions:
\begin{enumerate}
\item the curve $\Sigma$ is non-singular, i.e., all of its vertices  are of valence 3; 
\item for every $x\in H_1\cap\ldots\cap H_n$ which belongs to $\partial(C)\cap\Sigma$, the tuple $B_x=(\newton_x(H_1),\ldots, \newton_x(H_n))$ consists of a triangle and $n-1$ intervals.
\end{enumerate}
\end{rem}

\begin{defin}\label{interseps}
For an arbitrary $0\neq\varepsilon\in[-\varepsilon_1, \epsilon_1]$, by $\mathscr W(\varepsilon)$ we denote the finite set ${(C_1(\varepsilon)\cap H_2\cap\ldots\cap H_n)\setminus\{Q_1(\varepsilon),\ldots, Q_l(\varepsilon)\}}$. 
\end{defin}

Consider the sets $\mathscr W(-\varepsilon_1)=\{x_1(-\varepsilon_1),\ldots, x_p(-\varepsilon_1)\}$ and $\mathscr W(\varepsilon_1)=\{x_1(\varepsilon_1),\ldots, x_q(\varepsilon_1)\}$. 
We define 
\begin{equation*} 
\mathscr W(0)=\{\lim_{\varepsilon\to -0}(x_1(\varepsilon)), \ldots, \lim_{\varepsilon\to -0}(x_p(\varepsilon))\}\cup\{\lim_{\varepsilon\to +0}(x_1(\varepsilon)), \ldots, \lim_{\varepsilon\to +0}(x_q(\varepsilon))\}\subset (C\cap \Sigma).
\end{equation*}

\begin{rem}
It is straightforward to show that the sets $\mathscr W(0)$ and $\partial(C)\cap\Sigma$ coincide. 
\end{rem}

Condition 2 from Remark \ref{genericity} implies that for any point $x\in \mathscr W(0)=\partial(C)\cap\Sigma,$ there exist exactly two points $y_1=x_i(\pm\varepsilon_1), y_2=x_j(\pm\varepsilon_1)\in\mathscr W(-\varepsilon_1)\cup \mathscr W(\varepsilon_1)$ such that \\ $x=\lim_{\varepsilon\to\pm 0}(x_i(\varepsilon))=\lim_{\varepsilon\to\pm0}(x_j(\varepsilon))$. 
\begin{utver}
We have the following equality: 
\begin{multline}
{\mathrm{det}_2}(\newton_{y_1}(H_1(\pm\varepsilon_1)), \newton_{y_1}(H_2),\ldots\newton_{y_1}(H_n),\zeta)+\\{\mathrm{det}_2}(\newton_{y_2}(H_1(\pm\varepsilon_1)), \newton_{y_2}(H_2),\ldots\newton_{y_2}(H_n),\zeta)+\\{\mathrm{det}_2}(\newton_x(H_1),\ldots \newton_x(H_n),\zeta)=0.
\end{multline}
\end{utver}
$\vartriangleleft$ 
Consider the tuple $B_x=(\newton_x(H_1),\ldots, \newton_x(H_n))$. By the genericity assumption, it consists of a triangle and $n-1$ intervals. Moreover, $\mathsf K$ is the elementary obstacle corresponding to $B_x$. Since $x$ belongs to $\partial(C)\cap\Sigma$, by definition of the boundary of $C$, it follows that $\newton_x(H_i)$ is a triangle containing  $K_i$ as an edge for some $1\leqslant i\leqslant m$, while for every $j\neq i$ in $\{1,\ldots, m\}$ we have $\newton_x(H_j)=K_j$. Moreover, if $i=1$, then the other edges of the triangle are $\newton_{y_1}(H_1(\pm\varepsilon_1))$ and  $\newton_{y_2}(H_1(\pm\varepsilon_1))$. If $i\neq 1$, then these edges are $\newton_{y_1}(H_i)$ and $\newton_{y_2}(H_i)$.
Therefore, the sought equality follows from the linearity property of the 2-determinant. Indeed, by the genericity assumption, for $j\neq i$ the intervals $K_j,  \newton_{y_1}(H_j),  \newton_{y_2}(H_j)$ (or $K_1,  \newton_{y_1}(H_1(\pm\varepsilon_1)),  \newton_{y_1}(H_1(\pm\varepsilon_1))$, if $1=j\neq i$) coincide. Moreover, it is obvious that for $m+1\leqslant j\leqslant n$ the intervals $\newton_x(H_j), \newton_{y_1}(H_j)$ and $\newton_{y_2}(H_j)$ coincide, therefore, we have 
\begin{multline*}
{\mathrm{det}_2}(\newton_{y_1}(H_1(\pm\varepsilon_1)), \newton_{y_1}(H_2),\ldots\newton_{y_1}(H_n),\zeta)+\\{\mathrm{det}_2}(\newton_{y_2}(H_1(\pm\varepsilon_1)), \newton_{y_2}(H_2),\ldots\newton_{y_2}(H_n),\zeta)+{\mathrm{det}_2}(\newton_x(H_1),\ldots \newton_x(H_n),\zeta)=\\ {\mathrm{det}_2}(K_1,\ldots, K_{i-1}, \newton_{y_1}(H_i), K_{i+1}, \ldots, K_m, \newton_{y_1}(H_{m+1}), \ldots, \newton_{y_1}(H_n), \zeta)+\\{\mathrm{det}_2}(K_1,\ldots, K_{i-1}, \newton_{y_2}(H_i), K_{i+1}, \ldots, K_m, \newton_{y_2}(H_{m+1}), \ldots, \newton_{y_2}(H_n), \zeta)+\\{\mathrm{det}_2}(K_1,\ldots, K_{i-1}, K_i, K_{i+1}, \ldots, K_m, \newton_{x}(H_{m+1}), \ldots, \newton_{x}(H_n), \zeta)=0.
\end{multline*}
Exactly the same argument works in case $i=1$.
$\vartriangleright$

So, in order to prove Case 1 of Theorem \ref{calculS1}, it suffices to prove 

\begin{lemma}
The following equality holds:
$$\sum\limits_{x\in \mathscr W(0)}{\mathrm{det}_2}(K_1,\ldots, K_{i-1}, K_i, K_{i+1}, \ldots, K_m, \newton_{x}(H_{m+1}), \ldots, \newton_{x}(H_n), \zeta)=0.$$
\end{lemma}
$\vartriangleleft$ Each of the tuples $(K_1,\ldots, K_{i-1}, K_i, K_{i+1}, \ldots, K_m, \newton_{x}(H_{m+1}), \ldots, \newton_{x}(H_n), \zeta)$ has a subtuple $(K_1,\ldots, K_m)$ of rank $m-1$. 
There exists a natural projection $\pi\colon\Z^n\to\Z^n/\langle K_1,\ldots, K_{m}\rangle$. From Corollary \ref{blockdet2}, it follows that 
\begin{multline}
\sum\limits_{x\in \mathscr W(0)}{\mathrm{det}_2}(K_1,\ldots, K_{i-1}, K_i, K_{i+1}, \ldots, K_m, \newton_{x}(H_{m+1}), \ldots, \newton_{x}(H_n), \zeta)=\\ {\mathrm{det}_2}(K_1,\ldots, K_m)(\sum\limits_{x\in \mathscr W(0)}\det(\pi(\newton_{x}(H_{m+1})), \ldots, \pi(\newton_{x}(H_n)), \pi(\zeta))). 
\end{multline}
In order to show that the sum $\sum\limits_{x\in \mathscr W(0)}\det(\pi(\newton_{x}(H_{m+1})), \ldots, \pi(\newton_{x}(H_n)), \pi(\zeta))$ is equal to zero, we use Condition 1 from Remark \ref{genericity} and the assumption on the Newton Polytopes $\newton(H_1),\ldots, \newton(H_n)$ being 2-developed with respect to the point $\zeta$. The genericity assumption allows us to use the balancing condition at each of the vertices of the curve $\Sigma\cap C$. Fix an arbitrary smooth point $\mu$ on every edge of the curve $\Sigma\cap C$. For each of the vertices $a$ of this curve, consider the sum $M_a$ of the summands $\det(\pi(\newton_{\mu}(H_{m+1})), \ldots, \pi(\newton_{\mu}(H_n)), \pi(\zeta))$ over the edges adjacent to $a$. On the one hand, it immediately follows from the balancing condition that all such sums are equal to $0$. On the other hand, if we take the sum of $M_a$ over all the vertices of the curve $\Sigma\cap C$, we will obtain the sought sum. Indeed, since the summands over the edges with two endpoints in this sum are taken twice, they are cancelled. Thus, we obtain that the sought sum over the points $x\in\mathscr W(0)$ equals the sum over all rays (the edges with a single endpoint) of the classical intersection numbers of $C\cap\Sigma$ with the tropical hypersurface $S$ dual to the interval $E=\{0, \pi(\zeta)\}\subset\R^{n-m+1}$. This sum is always zero over $\mathbb F_2$. 

Indeed, pick any smooth point $x$ on a ray of the curve $\Sigma$, and denote the generating vector of this ray by $v$. 
By definition of a tuple of polytopes 2-developed with respect to a point (see Definition \ref{devwrt}), it follows that either $v(\zeta)=0$ over $\mathbb F_2$, or, for some $1\leqslant i\leqslant n$, the local Newton interval $\newton_x(\Sigma)$ has the endpoints of the same parity. In the first case, all the arguments of the sought determinant ${\mathrm{det}_2}(K_1,\ldots, K_{i-1}, K_i, K_{i+1}, \ldots, K_m, \newton_{x}(H_{m+1}), \ldots, \newton_{x}(H_n), \zeta)$ are contained in the same hyperplane $v^\perp$ in ${\mathbb F_2}^n$. In the second case, one of the arguments of the sougnh determinant is even, which finishes the proof of the lemma. 
$\vartriangleright$

Case 1 of Theorem \ref{calculS1} is proved. 

{\bf Case 2.} Consider a generic point $(H_1,\ldots, H_n)\in\EuScript{S}_1$ that belongs to a wall $L_{\mathsf K}$ corresponding to the elementary obstacle $\mathsf K=(K_i\mid i\in I)$ of Type 2. Without loss of generality suppose that $I=\{1,\ldots, m\}$ for some $m\leqslant n$. The elementary obstacle $K$ consists of a triangle and $m-1$ intervals. Moreover, we can assume that $K_1$ is the triangle. By $\Sigma$ denote the tropical curve obtained as the intersection $H_2\cap\ldots\cap H_n$. 

The set $H_1\cap\ldots\cap H_n$ consists of a finite set $\{Q_1,\ldots, Q_l\}$ of the transverse intersection points and the non-smooth point $x\in H_1$. Note that $x$ is a smooth point of the curve $\Sigma$. By $\overrightarrow{u}$ denote the primitive vector of the edge of $\Sigma$ which contains $x$. Pick any edge $K'_1$ of the triangle $K_1$ and consider the face $C_1$ of the hypersurface $H_1$ which is dual to $K'_1$. In analogy with the proof of Case 1, take any normal vector $\overrightarrow{w}$ to the hyperplane containing the face $C_1$. Moreover, we can assume that the vectors $\overrightarrow{u}$ and $\overrightarrow{w}$ are linearly independent, since we can always choose an edge $K'_1$ such that this condition is satisfied. 

Thus, we obtain the family of hypersurfaces  $H_1(\varepsilon)=T_{\varepsilon\overrightarrow{w}}(H_1)$ parametrized by $\varepsilon\in\R$. By $C_1(\varepsilon)$ we denote the image of the face $C_1$ of the hypersurface $H_1$.

The neighbourhood $V_0\ni (H_1,\ldots, H_n)$ constructed above contains an open ball $\EuScript B\ni (H_1,\ldots, H_n)$, therefore, there exists $\varepsilon_0>0$ such that for every $\varepsilon \in (-\varepsilon_0, \varepsilon_0)$ the point $(H_1(\varepsilon), H_2,\ldots, H_n)$ belongs to $\EuScript B$. Fix an arbitrary ${0<\varepsilon_1<\varepsilon_0}$. For every $\varepsilon\neq 0$ in $[-\varepsilon_1, \varepsilon_1]$, we have $H_1(\varepsilon)\pitchfork H_2\pitchfork\ldots\pitchfork H_n$. So, our next step is to compare the 2-intersection numbers $\iota_2(H_1(-\varepsilon_1), H_2, \ldots, H_n;\zeta)$ and $\iota_2(H_1(\varepsilon_1), H_2, \ldots, H_n;\zeta)$. 

From Proposition \ref{goodpoints}, it follows that for every transverse intersection point \\${Q_j\in H_1\cap H_2\cap\ldots\cap H_n}$, there exists a neighbourhood ${U_j\ni Q_j}$, such that for every ${0\neq\varepsilon\in [-\varepsilon_1, \varepsilon_1]}$, the neighbourhood $U_j$ contains a unique transverse intersection point ${Q_j(\varepsilon)\in H_1(\varepsilon)\cap H_2\cap\ldots\cap H_n}$. 

Besides the transverse intersection points $Q_j(\pm\varepsilon_1)$,  the union of the sets $\mathscr W(-\varepsilon_1)$ and $\mathscr W(\varepsilon_1)$ (see Definition \ref{interseps}) contains the smooth intersection points $y_1, y_2, y_3$ which appear as the result of the $\pm\varepsilon_1$-deformation of the non-smooth intersection \\${x\in H_1\cap\ldots\cap H_n}$. 

Therefore, for $\varepsilon=\pm\varepsilon_1$ we have the following equality:
\begin{multline}\label{iotasum2}
\iota_2(H_1(\varepsilon), H_2,\ldots, H_n;\zeta)= \\\sum_{j=1}^{j=l}{\mathrm{det}_2}(\newton_{Q_j(\varepsilon)}(H_1(\varepsilon)), \newton_{Q_j(\varepsilon)}(H_2),\ldots, \newton_{Q_j(\varepsilon)}(H_n),\zeta)+\\\sum\limits_{y\in \mathscr W(\varepsilon)}{\mathrm{det}_2}(\newton_y(H_1(\varepsilon)), H_2,\ldots \newton_y(H_n),\zeta). 
\end{multline}

By the construction of the interval $[-\varepsilon_1,\varepsilon_1]$, the tuples of local Newton intervals $B_{Q_j(-\varepsilon_1)}$ and $B_{Q_j(\varepsilon_1)}$ coincide for every $1\leqslant j\leqslant n$. Thus, to prove Case 2 of Theorem \ref{calculS1}, it suffices to prove the following 

\begin{lemma}
In the previous notation, the following equality holds: 
\begin{multline*}
\sum\limits_{y\in \mathscr W(-\varepsilon_1)}{\mathrm{det}_2}(\newton_y(H_1(-\varepsilon_1)), \newton_{y}(H_2),\ldots \newton_y(H_n),\zeta)+ \\\sum\limits_{y\in \mathscr W(\varepsilon_1)}{\mathrm{det}_2}(\newton_y(H_1(\varepsilon_1)), \newton_{y}(H_2),\ldots \newton_y(H_n),\zeta)=0.
\end{multline*}
\end{lemma}

$\vartriangleleft$ Rewrite the sought equality in terms of the intersection points $y_1,y_2, y_3$ defined above:
\begin{multline*}
{\mathrm{det}_2}(\newton_{y_1}H_1(\pm\varepsilon_1)), \newton_{y_1}(H_2),\ldots \newton_{y_1}(H_n),\zeta)+\\{\mathrm{det}_2}(\newton_{y_2}H_1(\pm\varepsilon_1)), \newton_{y_2}(H_2),\ldots \newton_{y_2}(H_n),\zeta)+\\{\mathrm{det}_2}(\newton_{y_3}H_1(\pm\varepsilon_1)),\newton_{y_3}(H_2),\ldots \newton_{y_3}(H_n),\zeta)=0
\end{multline*}
Note that the intervals $\newton_{y_1}H_1(\pm\varepsilon_1))$, $\newton_{y_2}H_1(\pm\varepsilon_1))$, $\newton_{y_3}H_1(\pm\varepsilon_1))$ are exactly the edges of the triangle $K_1$. Moreover, obviously, the other Newton intervals $\newton_{y_i}(H_j)$ coincide with $\newton_{x}(H_j)$ for all $i\in\{1,2,3\}$ and $j\in\{2,\ldots,n\}$. Therefore, the sought equality follows from the linearity property of the 2-determinant, which finishes the proof of the lemma.
$\vartriangleright$

So, Cases 1 and 2 of Theorem \ref{calculS1} are proved. The case of a tuple $H=(H_1,\ldots,H_n)$, corresponding to an elementary obstacle $K_1=\{a,b,c\}$ of type 3, is obvious: every intersection point of multiplicity ${\mathrm{det}_2}(a-c, \newton_{y_1}(H_2),\ldots \newton_{y_1}(H_n),\zeta)$ splits into two intersection points of multiplicites ${\mathrm{det}_2}(a-b, \newton_{y_1}(H_2),\ldots \newton_{y_1}(H_n),\zeta)$ and ${\mathrm{det}_2}(b-c, \newton_{y_1}(H_2),\ldots \newton_{y_1}(H_n),\zeta)$ as the tuple $H$ perturbs.
\end{proof}
Theorem \ref{calculS1} together with Proposition \ref{constS0} imply Theorem \ref{main}, therefore, we obtained a well-defined function ${\MV_2\colon(P_1,\ldots,P_n;\zeta)\mapsto\iota_2(H_1,\ldots, H_n;\zeta)}$ --- the so-called {\it 2-mixed volume}.

\section{Multivariate Vieta's Formula} \label{multVieta}
In this Section, we show that for a certain class of multivariate polynomial systems equations, the product of their roots can be expressed in terms of the 2-volume of their Newton polytopes. The Section is organised as follows. First we obtain such a formula for binomial systems (see Subsection \ref{BinomVieta}). In Subsection \ref{DefVieta}, we provide all the necessary definitions and notation crucial to formulating Theorem \ref{VietaMain} and conducting its proof. Subsection \ref{multvietaproof} is devoted to the statement and proof of the multivariate Vieta's formula (see Theorem \ref{VietaMain}).

\subsection{Multivariate Vieta's Formula for Binomial Systems}\label{BinomVieta}
\begin{lemma}\label{VietaBin}
Let $f_1,\ldots, f_n$ be binomials such that all the coefficients of $f_i, 1\leqslant i\leqslant n,$ are equal to 1 and the Newton intervals $\newton(f_1),\ldots,\newton(f_n)\subset\Z^n$ are affinely independent.Fix an arbitrary point $a\in \mathbb Z^n$. Then the following equality holds:
\begin{equation}
\prod_{f_1(x)=\ldots=f_n(x)=0}x^a=(-1)^{\MV_2(\newton(f_1),\ldots, \newton(f_n); a)} 
\end{equation}
\end{lemma}

$\vartriangleleft$ We will prove this lemma by induction on $n$. The base $n=1$ follows from the classical Vieta's formula for the product of roots for a polynomial in one variable. Now suppose that the statement is true for $k=n-1$. We will now deduce it in the case $k=n$.  

Since the product of roots for a system of equations is invariant under invertible monomial changes of variables and the Newton intervals of the polynomials of the system are affinely independent, we can assume without loss of generality, that the system $\{f_1=\ldots=f_n=0\}$ is of the following form:

\begin{equation}\label{binsys}
\begin{cases}
   {x_1}^{v_{1,1}}{x_2}^{v_{2,1}}\ldots {x_{n-1}}^{v_{{n-1},1}}+1=0\\
   {x_1}^{v_{1,2}}{x_2}^{v_{2,2}}\ldots {x_{n-1}}^{v_{{n-1},2}}+1=0\\
   \ldots \\
   {x_1}^{v_{1,{n-1}}}{x_2}^{v_{2,{n-1}}}\ldots {x_{n-1}}^{v_{{n-1},{n-1}}}+1=0\\
   {x_1}^{w_1} {x_2}^{w_2}\ldots{x_n}^{w_n}+1=0
 \end{cases}
\end{equation} 

Thus, the statement of Lemma \ref{VietaBin} for the system (\ref{binsys}) can be reformulated as follows: 
\begin{equation}
\prod_{f_1(x)=\ldots=f_n(x)=0}x^a=(-1)^{\det_2\left(\begin{smallmatrix}
v_{1,1} & \cdots & {v_{1,{n-1}}} & w_1 & a_1\\
v_{2,1} & \cdots & {v_{2,{n-1}}} & w_2 & a_2 \\         
\vdots & \ddots & \vdots & \vdots & \vdots \\
v_{{n-1},1} & \cdots & {v_{{n-1},{n-1}}} & w_{n-1} & a_{n-1} \\
0 & \cdots & 0 & w_n & a_{n}
\end{smallmatrix}\right)}\label{binVieta}
\end{equation}

By the induction hypothesis, we have the following equalities for $1\leqslant i\leqslant n-1$:

\begin{equation}\label{prodsmalli}
\prod_{f_1(x)=\ldots=f_{n-1}(x)=0} x_i=(-1)^{\det_2\left(\begin{smallmatrix}
v_{1,1} & \cdots & {v_{1,{n-1}}} & 0\\
\vdots & \ddots & \vdots & \vdots \\
v_{{i-1},1} & \cdots & {v_{{i-1},{n-1}}} & 0 \\
v_{i,1} & \cdots & {v_{i,{n-1}}} & 1 \\  
v_{{i+1},1} & \cdots & {v_{{i+1},{n-1}}} & 0 \\      
\vdots & \ddots & \vdots & \vdots \\
v_{{n-1},1} & \cdots & {v_{{n-1},{n-1}}} & 0 
\end{smallmatrix}\right)}
\end{equation} 

Using these equalities, it is easy to compute the product $\prod_{f_1(x)=\ldots=f_n(x)=0} x_i$. Indeed, each of the roots for the system $\{f_1(x)=\ldots=f_n(x)=0\}$ is obtained from substituting the roots for the smaller system $\{f_1(x)=\ldots=f_{n-1}(x)=0\}$ into the last equation $f_n(x)=0$ and solving this equation in the variable $x_n$. Thus, we have $w_n$ solutions for the system $\{f_1(x)=\ldots=f_n(x)=0\}$ corresponding to each of the roots $(\alpha_1, \ldots, \alpha_{n-1})$ for the system $\{f_1(x)=\ldots=f_{n-1}(x)=0\}$. Therefore, the sought product equals exactly the $w_n$-th power of the product (\ref{prodsmalli}). A straightforward computation consisting in applying the explicit formula for the 2-determinant (see Theorem \ref{formuladet2}) shows that the following equality holds for every $i\in\{1,\ldots, n-1\}$:
\begin{equation}
w_n {\mathrm{\det}_2\left(\begin{smallmatrix}
v_{1,1} & \cdots & {v_{1,{n-1}}} & 0\\
\vdots & \ddots & \vdots & \vdots \\
v_{{i-1},1} & \cdots & {v_{{i-1},{n-1}}} & 0 \\
v_{i,1} & \cdots & {v_{i,{n-1}}} & 1 \\  
v_{{i+1},1} & \cdots & {v_{{i+1},{n-1}}} & 0 \\      
\vdots & \ddots & \vdots & \vdots \\
v_{{n-1},1} & \cdots & {v_{{n-1},{n-1}}} & 0 
\end{smallmatrix}\right)}=\mathrm{\det}_2 \left(\begin{smallmatrix}
v_{1,1} & \cdots & {v_{1,{n-1}}} & w_1 & 0\\
\vdots & \ddots & \vdots & \vdots & \vdots \\
v_{{i-1},1} & \cdots & {v_{{i-1},{n-1}}} & w_{i-1} & 0 \\
v_{i,1} & \cdots & {v_{i,{n-1}}} & w_i & 1 \\  
v_{{i+1},1} & \cdots & {v_{{i+1},{n-1}}} & w_{i+1} & 0 \\      
\vdots & \ddots & \vdots & \vdots \vdots \\
v_{{n-1},1} & \cdots & {v_{{n-1},{n-1}}} & w_{n-1} & 0 \\
0 & \cdots & 0 & w_n & 0
\end{smallmatrix}\right)
\end{equation}

Now let us compute the product $\prod_{f_1(x)=\ldots =f_n(x)=0} x_n$. This can be easily done using the classical Vieta's formula. Consider a root $(\alpha_1,\ldots, \alpha_{n-1})$ for the system $\{f_1(x)=\ldots=f_{n-1}(x)=0\}$. Substituting it into the last equation $f_n(x)=0$ we obtain $w_n$ roots $(\alpha_1, \ldots, \alpha_{n-1}, \beta_j), 1\leqslant j\leqslant w_n,$  for the system $\{f_1(x)=\ldots =f_n(x)=0\}$. The product $\prod_{j=1}^{w_n} \beta_j$ equals $(-1)^{w_n}\alpha_1^{-w_1}\ldots\alpha_{n-1}^{-w_{n-1}}$. Therefore, we obtain that the sought product equals $(-1)^M$, where  
\begin{equation*}
M=w_n \vert\{x\mid f_1(x)=\ldots =f_{n-1}(x)=0\}\vert+{\prod_{f_1(x)=\ldots=f_{n-1}(x)=0} x_1^{-w_1}\ldots x_{n-1}^{-w_{n-1}}}.
\end{equation*}
From the Bernstein-Kushnirenko theorem (see \cite{bernstein}) and the induction hypothesis, we have 
\begin{equation}\label{binvietform}
M=w_n\det\left(\begin{smallmatrix}
v_{1,1} & \cdots & {v_{1,{n-1}}}\\
\vdots & \ddots & \vdots \\
v_{{n-1},1} & \cdots & {v_{{n-1},{n-1}}} 
\end{smallmatrix}\right)+\mathrm{\det}_2 \left(\begin{smallmatrix}
v_{1,1} & \cdots & {v_{1,{n-1}}} & w_1\\
\vdots & \ddots & \vdots & \vdots \\
v_{{n-1},1} & \cdots & {v_{{n-1},{n-1}}} & w_{n-1}\\
\end{smallmatrix}\right);
\end{equation}
one can easily check, that the expression in the right-hand side of the equality (\ref{binvietform}) equals exactly the following 2-determinant:
\begin{equation*}
M=\mathrm{\det}_2 \left(\begin{smallmatrix}
v_{1,1} & \cdots & {v_{1,{n-1}}} & w_1 & 0\\
\vdots & \ddots & \vdots & \vdots & \vdots \\
v_{{n-1},1} & \cdots & {v_{{n-1},{n-1}}} & w_{n-1} & 0 \\
0 & \cdots & 0 & w_n & 1
\end{smallmatrix}\right).
\end{equation*}
Lemma \ref{VietaBin} now follows, since the 2-determinant is multilinear and $\mathrm GL_n(\mathbb F_2)$-invariant.
$\vartriangleright$

The rest of the Section is devoted to a generalization of this result to a richer class of multivariate polynomial systems of equations.

\subsection{Some Necessary Notation and Definitions}\label{DefVieta}

Take an arbirtary point $0\neq a\in Z^n$ and consider an $a$--prickly tuple $P=(P_1,\ldots, P_n)$ of convex lattice polytopes in $\R^n$ (see Definition \ref{prickly}). By $\C_1^{P_i}$ we denote the set of all polynomials $f=\sum_{p\in P_i} c_p x^p$ such that $\newton(f)=P_i$ and if $p\in P_i$ is a vertex, then $c_p\ne 0$.
Consider the set $\C_1^{P}=\C_1^{P_1}\times\ldots\times\C_1^{P_n}$.  

\begin{defin}
Let $\gamma\neq 0$ in $(\R^*)^n$ be a covector. Let $f(x)$ be a Laurent polynomial with the Newton polytope $\newton(f)$. The {\it truncation} of $f(x)$ with respect to $\gamma$ is the polynomial $f^{\gamma}(x)$ that can be obtained from $f(x)$ by omitting the sum of monomials which are not contained in the support face 
${\newton(f)^{\gamma}}$. 
\end{defin}

It is easy to show that for a system of polynomial equations $\{f_1(x)=\ldots=f_n(x)=0\}$ and an arbitrary covector $\gamma\neq 0$, the ``truncated" system $\{f_1^{\gamma}(x)=\ldots=f_n^{\gamma}(x)=0\}$ by a monomial change of variables can be reduced to a system in $n-1$ variables at most. Therefore, for the systems with coefficients in general position, the ``truncated" systems are inconsistent in $\CC^n$. 

\begin{defin}
Let $(f_1,\ldots, f_n)$ be a tuple of polynomials in $\C_1^P$. In the same notation as above, the system ${f_1(x)=\ldots=f_n(x)=0}$ is called {\it degenerate at infinity}, if there exists a covector $\gamma\neq 0$ such that the system $\{f_1^{\gamma}(x)=\ldots=f_n^{\gamma}(x)=0\}$ is consistent in $\CC^n$.
\end{defin}

Let $\mathscr H\subset\C_1^P$ be the set of all systems that are degenerate at infinity. Consider $\mathscr H'\subset \mathscr H$ --- the set of the systems degenerate at infinity, which satisfy the following property: there exists a covector $\gamma\in(\R^*)^n$ such that the system $\{f_1^{\gamma}(x)=\ldots=f_n^{\gamma}(x)=0\}$ is consistent in $(\C^*)^n$ and the Minkowski sum of the support faces ${\newton(f_i)^{\gamma}}, 1\leqslant i\leqslant n$ is of codimension greater than~$1$. 
One can easily see that $\codim(\mathscr H')>1$. 

We denote by $\mathscr D$ the discriminant hypersurface in $\C_1^P$ (i.e. the closure of the set of all systems that have a multiple root in $\CC^n$), and $\mathscr D'\subset\mathscr D$ stands for the set of all $(f_1,\ldots f_n)\in\C_1^P$ such that the system $\{f_1(x)=\ldots=f_n(x)=0\}$ has non-isolated roots. It is easy to show that $\codim(\mathscr D')>1$.

\begin{defin}
The system $\{f_1(x)=\ldots=f_n(x)=0\}$ with $F=(f_1,\ldots, f_n)\in\C_1^P$ is called {\it non-degenerate}, if $F\in\C_1^P\setminus(\mathscr H\cup\mathscr D)$. 
\end{defin}

\subsection{The Multivariate Vieta's Formula: Statement and Proof}\label{multvietaproof}

Fix an arbitrary point $0\neq a\in\Z^n$ and an $a$-prickly tuple $P=(P_1,\ldots, P_n)$ of polytopes. In the notation of Subsection \ref{DefVieta}, the multivariate Vieta's formula expresses the product of the monomials $x^a$ over all the roots $x$ for a system of polynomial equations $f_1(x)=\ldots=f_n(x)$, where $F=(f_1,\ldots,f_n)\in\C_1^P$ and the coefficients of $f_i$ at the vertices of its Newton polytope are equal to 1, in terms of the 2-mixed volume (see Section \ref{2-volume}) of the polytopes $P_1,\ldots, P_n$ and the point $a$. This Subsection consists of two parts. First, we obtain a holomorphic everywhere defined function $\Phi\colon\C_1^P\to\C$ with no zeroes and poles, which maps a point $(f_1,\ldots,f_n)\in\C_1^P$ to the abovementioned product of monomials. It immediately follows that the function $\Phi$ is constant on tuples $(f_1,\ldots,f_n)$ such that the coefficients of $f_i$ at the vertices of its Newton polytope are equal to 1. Secondly, we compute this constant, reducing this problem to the case of binomial systems of equations, which was studied in Subsection \ref{BinomVieta}.

In the notation of Subsection \ref{DefVieta}, consider the function $\Phi_0\colon\C_1^P\setminus(\mathscr H\cup\mathscr D)\to\C$, defined as follows: 
\begin{equation*}
\Phi_0\colon (f_1,\ldots, f_n)\mapsto\prod_{f_1(x)=\ldots=f_n(x), x\in\CC^n} x^a.
\end{equation*}

\begin{theor}\label{const}
The function $\Phi_0$ is a monomial in the coefficients of $f_1,\ldots,f_n$ at the vertices of their Newton polytopes. 
\end{theor}

\begin{proof}
Obviously, $\Phi_0$ is a well-defined holomorphic function with no zeroes and poles in the open dense subset $\C_1^P\setminus(\mathscr H\cup\mathscr D)$. Our goal is to show that this function can be extended to a holomorphic function $\Phi\colon\C_1^P\to\C$ with no zeroes and poles and, therefore, is a monomial in the coefficients of $f_1,\ldots,f_n$ at the vertices of their Newton polytopes. We will construct this extension in several steps.  
\begin{lemma}
There exists a holomorphic extension of the function $\Phi_0$ on the set ${\C_1^P\setminus(\mathscr H'\cup\mathscr D)}$.
\end{lemma}
$\vartriangleleft$ 
Since the product of roots is invariant under monomial changes of variables, we can assume without loss of generality that $a=(0,\ldots,0,\alpha)\in\Z^n$. Consider the projection $\pi\colon\Z^n\to\Z^{n-1}$ along the radius vector of the point $a$ and the tuple $Q=(Q_1,\ldots,Q_n)$ of convex lattice polytopes, where $Q_i=\conv(\pi(P_i))\subset\R^{n-1}$. Take a simple fan $\Omega$ compatible with the tuple $Q$ of polytopes and consider the smooth projective toric variety $X_{\Omega}$. We obtain an inclusion $\CC^{n-1}\times\CC\hookrightarrow X_{\Omega}\times\CC$. 

Let $\Phi_1\colon\C_1^p\setminus(\mathscr H'\cup\mathscr D)\to\C$ be the function defined as follows:
\begin{equation}\label{rootprod}
\Phi_1\colon (f_1,\ldots, f_n)\mapsto \prod_{f_1(x)=\ldots=f_n(x), x\in X_{\Omega}\times\CC} x^a.
\end{equation}
This function is well-defined, since at a point $x=(y,t)\in X_{\Omega}\times\CC$, the monomial $x^a$ equals $t^{\alpha}\in\CC$, therefore, the function $\Phi_1$ is the sought extension of $\Phi_0$. $\vartriangleright$
\begin{lemma}
The function $\Phi_1$ can be regularly extended to a function on the set ${\C_1^p}$.
\end{lemma} 
$\vartriangleleft$ 
\begin{utver}\label{killsing}
Let $\Psi\colon\C^N\setminus\Sigma\to\mathbb C$, where~ $\codim (\Sigma)=1$, be a holomorphic function. If $\Sigma'\subset \Sigma$~is such that $\codim(\Sigma')>1$ and  $\Psi$ can be continuously extended to $\tilde{\Psi}\colon\C^N\setminus\Sigma'\to\C$, then there exists a holomorphic extension $\bar{\Psi}\colon\mathbb{C}^N\to\mathbb C$ of $\Psi$.
\end{utver}

Take an isolated multiple root $q$ for the system $\{f_1=\ldots=f_n=0\}$, where $(f_1,\ldots, f_n)\in\C_1^P\setminus(\mathscr H'\cup \mathscr D')$. Its multiplicity equals the degree of the map $F\colon\mathscr U\to\C^n$, where $\mathscr U\subset \C^n$ is a small neighborhood of the point $q$, $F\colon x\mapsto (f_1(x),\ldots, f_n(x))$. For almost all points $\varepsilon=(\varepsilon_1,\ldots, \varepsilon_n)$ in a sufficiently small neighborhood $0\in\mathscr V\subset\C^n$, the number of preimages $F^{-1}(\varepsilon)$ equals the same number $k$ of multiplicity $1$ roots $q_1(\varepsilon),\ldots, q_k(\varepsilon)$ of the system $\{f_1(x)=\varepsilon_1,\ldots f_n(x)=\varepsilon_n\}$, which are contained in $\mathscr U$. Moreover, as $\varepsilon_i\rightarrow 0, 1\leqslant i\leqslant n$, we have $q_j(\varepsilon)\rightarrow q, 1\leqslant j\leqslant k$. Therefore, for the product of monomials, we have $\prod_1^k (q_k(\varepsilon))^a\rightarrow (q^k)^a$. Thus, letting the monomials $x^a$ enter the product the number of times equal to the multiplicity of the corresponding root $x$, we obtain a continuous extension $\tilde{\Phi}_1\colon\C_1^P\setminus (\mathscr H'\cup\mathscr D')$. 

Using Proposition \ref{killsing} for $\Psi=\Phi_1$, $\Sigma=\mathscr H\cup\mathscr D$ and $\Sigma'=\mathscr H'\cup\mathscr D'$, we obtain the desired holomorphic extension $\Phi\colon\C_1^P\to\C$. $\vartriangleright$

Note that the function $\Phi$ has no zeroes and poles in $\C_1^P$, therefore, it is a monomial, which finishes the proof of the theorem. 
\end{proof}

It follows from Theorem \ref{const} that there exists a well-defined function $\Phi(P_1,\ldots, P_n; a)$, which maps $a\in\Z^n$ is a point and $P=(P_1,\ldots, P_n)$ is an $a$-prickly tuple of polytopes in $\R^n$ to the product (\ref{rootprod}) of monomials over the roots for a polynomial system of equations ${\{f_1(x)=\ldots=f_n(x)=0\},}$ where ${(f_1,\ldots, f_n)\in\C_1^P}$ is an arbitrary point such that all coefficients of $f_i$ at the vertices of its Newton polytope are equal to 1. 

Now we are ready to state the multivariate Vieta's formula.
\begin{theor}\label{VietaMain}
Under the same assumptions as above, we have
\begin{equation}
\Phi(P_1,\ldots, P_n; a)=(-1)^{\MV_2(P_1,\ldots, P_n; a)}. \label{Vieta}
\end{equation}
\end{theor}

The main idea of the proof is to introduce a new variable, a parameter $t$, in such a way that as $t\rightarrow\infty$, the system of equations asymptotically breaks down into a union of binomial systems, which were considered in Subsection \ref{BinomVieta}. 

Consider an arbitrary point $(f_1,\ldots, f_n)\in\C_1^P$, where $f_i(x)=\sum_{k\in P_i\cap\Z^n} c_{i,k}x^k$. With the system $\{f_1(x)=\ldots=f_n(x)=0\}$ we associate a perturbed system $\{\tilde{f}_1(x,t)=\ldots=\tilde{f}_n(x,t)=0\}$ of equations of the following form: $$\tilde{f}_i(x,t)=\sum_{k\in P_i\cap\Z^n} c_{i,k}t^{n_{i,k}}x^k, $$ where $n_{i,k}$ are non-negative integers. 
The Newton polytopes of $\tilde{f}_i, 1\leqslant i\leqslant n,$ are denoted by $\tilde{P}_i$. Note that $\tilde{P}_i\subset\R^n\times \R$ are lattice polytopes lying over the polytopes $P_i$. By $\rho\colon\R^{n}\times\R\to\R^n$ we denote the projection forgetting the last coordinate. 

\begin{defin}\label{uppface}
A face $\tilde{\Gamma}\subset\tilde{P}_i$ is said to be an {\it upper face}, if there exists a covector $v\in(\R^*)^n+1$ with the strictly positive last coordinate such that $\tilde{\Gamma}=\tilde{P}_i^v$.
\end{defin}

Each of the polytopes $\tilde{P}_i$ defines a convex subdivision $\Delta_i$ of the polytope $P_i$. Namely, each cell $\Gamma\in\Delta_i$ is equal to $\rho(\tilde{\Gamma})$ for some upper face $\tilde{\Gamma}\subset\tilde{P}_i$. In the same way, the Minkowski sum $\EuScript{\tilde{P}}=\sum \tilde{P}_i$ defines a convex subdivision $\Delta$ of the polytope $\EuScript{P}=\sum P_i$. Each cell $\Gamma\in\Delta$ can be uniquely represented as a sum $\Gamma=\Gamma_1+\ldots+\Gamma_n$, where $\Gamma_i\in\Delta_i$. 

To the tuple of polynomials $\tilde{f}_1,\ldots, \tilde{f}_n, \tilde{f}_i=\sum_{k\in P_i\cap\Z^n} c_k x^k t^{n_{i,k}}$, we associate tuples of polynomials $g_1,\ldots, g_n$ and $\tilde{g}_1,\ldots, \tilde{g}_n$ defined as follows: $$g_i=\sum_{k\in P_i\cap\Z^n} a_k x^k,$$ where $a_k$ equals $1$, if $k$ is a vertex of a cell of $\Delta_i$, and $0$ otherwise, and $\tilde{g}_i=\sum a_k x^k t^n_k, 1\leqslant i\leqslant n$. 
From Theorem \ref{const}, it follows that the sought product of monomials does not depend on the choice of the system of equations, so, it suffices to compute it for the system $\{g_1(x)=\ldots=g_n(x)=0\}$.

\begin{defin}\label{consistface}
A tuple of faces $\Gamma_1,\ldots,\Gamma_n, \Gamma_i\subset\tilde{P}_i,$ is said to be {\it consistent}, if there exists a covector $v\in(\R^*)^n$ such that $\Gamma_i=\tilde{P}_i^{v}$.
\end{defin} 

\begin{defin}\cite{det2}\label{defaffin}
A tuple of consistent faces $\Gamma_1,\ldots,\Gamma_n,  \Gamma_i\subset\tilde{P}_i,$ is said to be {\it affinely independent}, if for the Minkowski sum $\Gamma=\Gamma_1+\ldots+\Gamma_n$, the equality $\dim(\Gamma)=\dim(\Gamma_1)+\ldots+\dim(\Gamma_n)$ holds. 
\end{defin}
\begin{rem}
Clearly, if the faces $\Gamma_1,\ldots, \Gamma_n$ are affinely independent, then either one of these faces is a point, or these faces are all segments.
\end{rem}
\begin{utver}\cite{det2}\label{affind}
In the previous notation, we can choose the polytopes ${\tilde{P_i},1\leqslant i\leqslant n,}$ in such a way that any consistent tuple $\Gamma_1,\ldots, \Gamma_n$ of faces is affinely independent. 
\end{utver}

It follows from Proposition \ref{affind} that without loss of generality we can assume that all the consistent tuples $\Gamma_1,\ldots, \Gamma_n$ are affinely independent. 

By $\K$ we denote the field $\C\{\{t\}\}$ of Puiseux series with the standard valuation function $\val\colon(\C\{\{t\}\}\setminus 0)\to\R$. If $0\neq b\in \K$ and $\val(b)=q$, then we write $b=\beta t^q+\ldots$ to distinguish the leading term. Note also that we identify the space $\R^n$ with its dual by means of the standard inner product.  

For almost all the values $\tau$ of the parameter $t$, the systems ${\{\tilde{g}_1(x,\tau)=\ldots=\tilde{g}_n(x,\tau)=0\}}$ have the same finite number of roots $z_1(\tau),\ldots, z_d(\tau)$. It follows from the Bernstein--Kushnirenko formula that $d=\MV(P_1,\ldots, P_n)$, see \cite{bernstein} for the details. Each of the roots $z_j(t)$ can be considered as an algebraic function $z_j\colon\CC\to\CC^n$ in the variable $t$. The Puiseux series of the curve $z_j$ at infinity is of the following form:
\begin{equation}\label{roott}
z_j(t)=(\beta_1t^{v_1},\ldots,\beta_n t^{v_n})+\mathrm{componentwise~lower~ terms}. 
\end{equation}
For simplicity of notation, we shall write the expression (\ref{roott}) in the following form: $z_j(t)=(B_j t^{V_j})+\ldots$, where $B_j=(\beta_1,\ldots, \beta_n)$ and $V_j$ is the valuation vector $(v_1,\ldots, v_n)$. 

In the previous notation, consider the valuation vectors $V_j=(v_1\ldots, v_n)$ and $W_j=(v_1, \ldots, v_n, 1)$ of the roots $z_j(t)$ and $Z_j=(z_j(t), t)$ respectively. The function $\langle W_j, \cdot \rangle$ attains its maximum at some faces $\Gamma_1^{W_j}, \ldots,\Gamma_n^{W_j}, \Gamma_i^{W_j}\subset \tilde{P}_i$, which are, according to our assumption, affinely independent, and, therefore, if none of them is a vertex, then they are all segments. By $G_i^{W_j}$ we denote the set $\rho(\Gamma_i^{W_j})\cap \Z^n$. Note that if $\Gamma_i^{W_j}$ is a segment, then $\vert G_i^{W_j} \vert=2$ (those 2 elements are exactly the endpoints of the segment $\rho(\Gamma_i^{W_j})$).

Now, let us substitute the root $Z_j=(z_j(t),t)$ into the equations $\tilde{g}_1(x,t), \ldots, \tilde{g}_1(x,t)$. On one hand, what we obtain is nothing but $0$. On the other hand, the result is a Puiseux series in the variable $t$:
\begin{equation}\label{rootpuiseux}
\tilde{g}_i(z_j(t),t)=\sum_{k\in P_i\cap\Z^n} a_k (\prod_{m=1}^n \beta_m^{v_m}) t^{\langle V_j, k \rangle + n_{i,k}}+\ldots=\sum_{k\in P_i\cap\Z^n} a_k B_j^{V_j} t^{\langle V_j, k \rangle + n_{i,k}}+\ldots.
\end{equation}

\begin{rem}\label{binreduce}
It is obvious that the leading coefficient of the series (\ref{rootpuiseux}) equals the sum $$h_i(B_j)=\sum_{k\in G_i^w} a_k B_j^k,$$ while $\tilde{g}_i(z_j(t),t)=0$. Therefore, $B_j=(\beta_1, \ldots, \beta_n)$ is a root of the system $\{h_1(x)=\ldots=h_n(x)=0\}$, where each of the polynomials $h_i$ is obtained by omitting the terms that are not contained in $G_i^{W_j}$. 
\end{rem}

\begin{rem}
If the tuple $G_i^{W_j}, 1\leqslant j\leqslant n,$ contains a singleton $G_m^{W_j}$ for some $m$, then, the corresponding system $\{h_1(x)=\ldots=h_n(x)=0\}$ is inconsisent in $\CC^n$, since the polynomial $h_m$ is a monomial. Therefore, all the roots $B_j$ mentioned in Remark \ref{binreduce} are the roots for the corresponding binomial system $\{h_1(x)=\ldots=h_n(x)=0\}$. The number of the roots for such a system equals $\MV(G_1^{W_j}, \ldots, G_n^{W_j})$, by the Bernstein--Kushnirenko theorem, see \cite{bernstein}. 
\end{rem}

\begin{proof}[Proof of Theorem \ref{VietaMain}.]
It follows from Remark \ref{binreduce} and Lemma \ref{VietaBin} that the product of monomials from Theorem \ref{VietaMain} equals the sign of the limit 
\begin{equation*}
\lim_{t\rightarrow+\infty}(\prod_{j=1}^d(Z_j(t))^a)=\lim_{t\rightarrow+\infty}(\prod_{j=1}^d (B_j t^{V_j})^a+\ldots)=\lim_{t\rightarrow+\infty}(\prod_{j=1}^d B_j^a t^{\langle V_j,a\rangle})+\ldots).
\end{equation*}

Note that by Lemma \ref{VietaBin} and the definition of the 2-mixed volume (see Definition \ref{2-volume}), the product $\prod_{j=1}^d(B_j)^a$ equals exactly the 2-mixed volume $(-1)^{\MV_2(P_1,\ldots, P_n; a)}$, so, the sought product can be expressed as the limit 
\begin{equation}\label{rootmove}
\lim_{t\rightarrow+\infty}(\prod_{j=1}^d(Z_j(t))^a)=(-1)^{\MV_2(P_1,\ldots, P_n; a)}\lim_{t\rightarrow+\infty}(t^{\sum^{d}_{j=1} \langle V_j,a\rangle}+\ldots).
\end{equation}

From Theorem \ref{const}, it follows that the sought product of monomials is a monomial.  
Therefore, it suffices to compute the sign of the limit (\ref{rootmove}) as $t\rightarrow\infty, t\in \R_{>0}$, which  obviously equals $(-1)^{\MV_2(P_1,\ldots, P_n; a)}$. 
\end{proof}

\section{Signs of the Leading Coefficients of the Resultant} \label{resultant}

In this Section, we show how to reduce the computation of the leading coefficients of the resultant to finding the product of roots for a certain system of equations. In Subsection \ref{NewtPolRes}, we recall the definition of the {\it sparse mixed resultant}. Then, using the multivariate Vieta's formula (see Theorem \ref{VietaMain}), we obtain a closed-form expression in terms of the 2-mixed volume to compute the signs of the leading coefficients of the sparse mixed resultant. 

\subsection{The Sparse Mixed Resultant and Its Newton Polytope}\label{NewtPolRes}
Recall the definition of the codimension of a tuple of finite sets in $\Z^n$ given in Subsection \ref{buildwall}. 

\begin{defin} 
Let $A=(A_0,\ldots, A_n)$ be a $n$-tuple of finite sets in $\Z^n$. For an arbitrary subset $I\subset\{0,\ldots,n\}$, we define its {\it codimension}, which we denote by $\codim(I)$, as follows:
$$ \codim(I)=\dim(\sum_{i\in I} A_i)-|I|. $$
The {\it codimension of the tuple} $A$ is defined by the following equality: 
$$\codim(A)=\min_{I\subset\{0,\ldots,n\}}\codim(I).$$
\end{defin}

\begin{defin}\cite{sturmfels}\label{resuldef}
Consider a tuple $A=(A_0,\ldots, A_n)$ of finite sets in $\Z^n$ such that $\codim(A)=-1$ and the sets $A_i$ jointly generate the affine lattice $\Z^n$. Then the {\it sparse mixed resultant} $\resul(A)$ is a unique (up to scaling) irreducible polynomial in $\sum_0^n|A_i|$ variables $c_{i,a}$ which vanishes whenever the Laurent polynomials $f_i(x)=\sum_{a\in A_i} c_{i,a}x^a$ have a common zero in $\CC^n$.
\end{defin}

Here we provide an explicit description of the vertices of the Newton polytope of the sparse mixed resultant $\resul(A)$. For more details and proofs we refer the reader to the paper \cite{sturmfels}.

Let $Q=(Q_0,\ldots, Q_n)$ be the tuple of convex hulls of the sets $A_i$. Let $\omega\colon\cup_{i=0}^n Q_i\to \R_{\geqslant 0}$ be an arbitrary function. By $\tilde{Q}=(\tilde{Q}_0,\ldots, \tilde{Q}_n)$ denote the tuple of convex hulls of the sets $Q_i(\omega)=\{(a, \omega(a))\mid a\in Q_i\}$ and by $\rho$ the standard projection $\R^{n+1}\to\R^n$ forgetting the last coordinate. The polytopes $\tilde{Q}_i, 0\leqslant i\leqslant n,$ define convex subdivisions $\Delta_i$ and $\Delta$ of polytopes $Q_i, 0\leqslant i\leqslant n,$ and $\EuScript{Q}=Q_0+Q_1+\ldots+Q_n$, respectively. Namely, each cell $\Gamma\in\Delta_i$ is equal to $\rho(\tilde{\Gamma})$ for some upper face $\tilde{\Gamma}\subset\tilde{Q}_i$. In the same way, the Minkowski sum  $\EuScript{\tilde{Q}}=\sum \tilde{Q}_i$ defines a convex subdivision $\Delta$ of the polytope $\EuScript{Q}=\sum Q_i$. 

\begin{defin}
In the previous notation, a tuple of subdivisions $\Delta_i$ of the polytopes $Q_i$ is called a {\it mixed subdivision} of the polytopes $Q_0,\ldots, Q_n$. 
\end{defin}

Note that each upper face $\tilde{\Gamma}\subset\EuScript{\tilde{Q}}$ (see Definition \ref{uppface}) can be uniquely represented as a sum ${\tilde{\Gamma}=\tilde{\Gamma}_0+\ldots+\tilde{ 
\Gamma}_n}$, where $\tilde{ \Gamma}_i\subset\tilde{Q}_i$. Proposition \ref{affind} implies that without loss of genericity, we can assume that the all the tuples of consistent faces $(\tilde{\Gamma}_0,\ldots, \tilde{\Gamma}_n)$ are affinely independent (see Definitions \ref{consistface} and \ref{defaffin}). By the Dirichlet principle, we have that for every upper face ${\tilde{\Gamma}}\subset \EuScript{\tilde{Q}}$, one of the faces, say, $\tilde{ \Gamma}_j$, employed in the decomposition  ${\tilde{\Gamma}=\tilde{\Gamma}_0+\ldots+\tilde{
\Gamma}_n}$ is a vertex. By $c_{\Gamma}$ we denote the corresponding coefficient $c_{j, \rho{\tilde{\Gamma}}}$ of the polynomial $f_j$ of the system $\{f_0(x)=\ldots=f_n(x)=0\}$. 

To each of the mixed subdivisions we can associate a vertex of the Newton polytope $\newton{\resul(A)}$ as follows. The leading term of $\resul(A)$ which corresponds to a given mixed subdivision $\Delta_0,\ldots,\Delta_n$ is the product over all the facets $\Gamma\in \Delta$ (i.e., faces of codimension 1) in the subdivision $\Delta$ of the Minkowski sum $\EuScript{Q}$ of the multiples equal to $c_{\Gamma}^{\Vol(\Gamma)}$. 

\begin{theor}[B.Sturmfels, \cite{sturmfels}]
The construction given above provides a bijection between the set of all the mixed subdivisions of the tuple $Q=(Q_0,\ldots, Q_n)$ and the vertices of the Newton polytope $\newton(\resul(A))$ of the resultant.
\end{theor}

\begin{exa}
Take $A=(A_0, A_1)$, where $A_0=\{0,1\}, A_1=\{0,1,2\}\subset\Z$. Our aim is construct the Newton polytope of the polynomial $\resul(A)\in\C[a_0, a_1, b_0, b_1, b_2]$ which vanishes whenever the system $\{a_0+a_1x=b_0+b_1x+b_2x^2=0\}$ is consistent. In the same notation as above, we have $Q_0=[0,1], Q_1=[0,2]$. The following figure describes one of the three possible mixed subdivisions of the intervals $Q_0$ and $Q_1$ which yields the vertex $(1,1,0,1,0)\in\newton(\resul(A))\subset\R^5$.

\begin{center}
\begin{tikzpicture}
\draw[help lines] (0,0) grid (2,2);
\draw[help lines] (3,0) grid (5,2);
\draw[help lines] (7,0) grid (11,4);
\draw[ultra thick, violet] (7,0)--(8,0);
\draw[ultra thick, orange] (8,0)--(9,0);
\draw[ultra thick, green] (9,0)--(10,0);
\draw[fill] (3,0) circle[radius=0.07];
\draw[fill] (5,0) circle[radius=0.07];
\draw[fill] (7,0) circle[radius=0.07];
\draw[fill] (8,0) circle[radius=0.07];
\draw[fill] (9,0) circle[radius=0.07];
\draw[fill] (10,0) circle[radius=0.07];
\draw[ultra thick] (0,0)--(0,1);
\draw[ultra thick] (1,0)--(1,1);
\draw[ultra thick] (3,0)--(3,1);
\draw[ultra thick] (4,0)--(4,2);
\draw[ultra thick] (5,0)--(5,1);
\draw[ultra thick] (7,0)--(7,2);
\draw[ultra thick] (8,0)--(8,3);
\draw[ultra thick] (9,0)--(9,3);
\draw[ultra thick] (10,0)--(10,2);
\draw[ultra thick, violet] (7,2)--(8,3);
\draw[ultra thick, orange] (8,3)--(9,3);
\draw[ultra thick, green] (9,3)--(10,2);
\draw[fill] (7,2) circle[radius=0.07];
\draw[fill] (8,3) circle[radius=0.07];
\draw[fill] (9,3) circle[radius=0.07];
\draw[fill] (10,2) circle[radius=0.07];
\draw[ultra thick, violet] (3,1)--(4,2);
\draw[ultra thick, orange] (0,1)--(1,1);
\draw[ultra thick, green] (4,2)--(5,1);
\draw[thick] (3,0)--(4,0);
\draw[thick] (0,0)--(1,0);
\draw[thick] (4,0)--(5,0);
\node[below] at (0,-0.1) {$a_0$};
\node[below] at (1,-0.1) {$a_1$};
\node[below] at (3,-0.1) {$b_0$};
\node[below] at (4,-0.1) {$b_1$};
\node[below] at (5,-0.1) {$b_2$};
\draw[fill, violet] (0,1) circle[radius=0.1];
\draw[fill, orange] (4,2) circle[radius=0.1];
\draw[fill, green] (1,1) circle[radius=0.1];
\draw[fill] (3,1) circle[radius=0.07];
\draw[fill] (5,1) circle[radius=0.07];
\draw[fill, violet] (0,0) circle[radius=0.14];
\draw[fill, green] (1,0) circle[radius=0.14];
\draw[fill, orange] (4,0) circle[radius=0.14];
\draw [|->, ultra thick, red] (12,1.5) -- (13.5,1.5);
\node [right, violet] at (14,1.5) {$a_0^1$};
\node [right, green] at (14.4,1.5) {$a_1^1$};
\node [right, orange] at (14.8,1.5) {$b_1^1$};
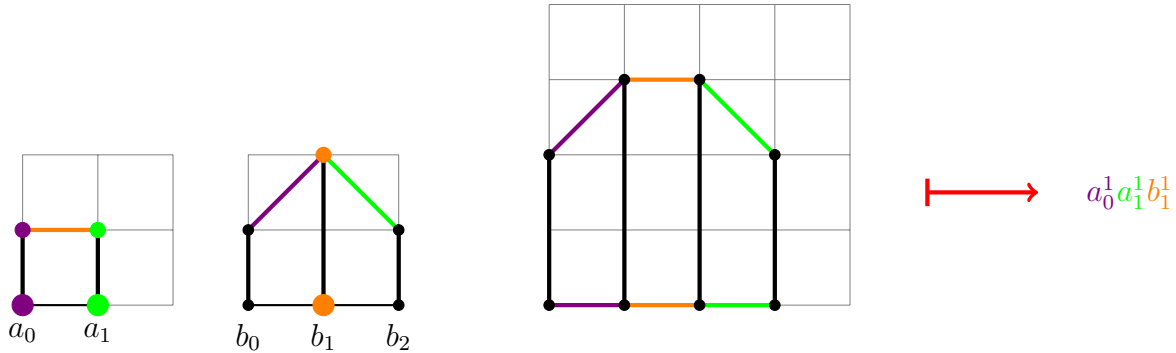
\node[below] at (7, -1.2) {Figure 3. The bijection between the mixed subdivisions of $A$ and the vertices of $\newton(\resul(A))$.};
\end{tikzpicture} 
\end{center}
\end{exa}

\subsection{Computing the Signs of the Leading Coefficients of the Resultant}\label{SignsRes}

In the previous notation, consider a tuple $A=(A_0,\ldots, A_n)$ of finite sets in $\Z^n$ satisfying the properties given in Definition \ref{resuldef}. Recall that by $|A_i|$ we denote the cardinality of the set $A_i\subset \Z^n$, and $|A|$ stands for the sum $\sum_0^n|A_i|$. For simplicity of notation, by $\resul$ we  denote the sparse mixed resultant $\resul(A)$.  

Consider the Newton polytope $\newton(\resul)$ of the resultant $\resul(A)$ (see \cite{sturmfels} for its explicit description). Suppose that we are given a pair of gradings ${\gamma=(\alpha_{i,a}\mid i\in\{0,\ldots,n\}, a\in A_i)}$ and ${\sigma=(\beta_{j,b}\mid j\in\{0,\ldots,n\}, b\in A_i)\in (\Z^*)^{|A|}}$ with strictly positive coordinates such that the support faces $\newton(\resul)^{\gamma}$ and $\newton(\resul)^{\sigma}$ are $0$-dimensional. We will now compute the quotient of the coefficients ${r_{\gamma}}$ and ${r_{\sigma}}$ of $\resul$ which are leading with respect to the gradings $\gamma$ and $\sigma$ respectively, by reducing this problem to the multivariate Vieta's formula (see Theorem \ref{VietaMain}).

To the covectors $\gamma, \sigma$ one can associate the tuple $P^{\gamma,\sigma}=(P_0^{\gamma,\sigma},\ldots, P_n^{\gamma,\sigma})$ of polytopes in $\R^{n+1}$ such that \\ $${P_i^{\gamma,\sigma}=\conv(\{(a, \alpha_{i,a})\mid a\in A_i\}\cup\{(a, -\beta_{i,a})\mid a\in A_i\})}.$$

\begin{exa}\label{exa1}
Let $A=(A_0, A_1)$, where $A_0=\{0,1\}, A_1=\{0,1,2\}\subset\Z$. The Newton polytope $\newton(\resul(A))$ is a triangle with vertices $\bar{\gamma}=(2,0,0,0,1)$, $\bar{\sigma}=(0,2,1,0,0)$ and $\bar{\delta}=(1,1,0,1,0)$. Consider the covectors $\gamma=(2,1,1,1,2)$, $\sigma=(1,2,2,1,1)$, and $\delta=(2,2,1,2,1)$, whose support faces are the vertices $\bar{\gamma}, \bar{\sigma}, \bar{\delta}$. Thus, we obtain the the polygons $P_0^{\gamma,\sigma}$ and $P_1^{\gamma,\sigma}$ (see fig. 4a) and the polygons $P_0^{\gamma,\delta}$ and $P_1^{\gamma,\delta}$ (see fig. 4b).
\end{exa}

\begin{center}
\begin{tikzpicture}
\draw[ultra thick](-5,0)--(-4,0);
\draw[fill] (-5,0) circle[radius=0.06];
\draw[fill] (-4,0) circle[radius=0.06];
\draw[ultra thick] (-5,-1)--(-5,2);
\draw[fill] (-5,-1) circle[radius=0.1];
\draw[fill] (-5,2) circle[radius=0.1];
\draw[ultra thick] (-4,-2)--(-4,1);
\draw[fill] (-4,-2) circle[radius=0.1];
\draw[fill] (-4,1) circle[radius=0.1];
\draw[ultra thick] (-5,2)--(-4,1);
\draw[ultra thick] (-5,-1)--(-4,-2);
\draw[ultra thick](-1,0)--(1,0);
\draw[fill] (-1,0) circle[radius=0.06];
\draw[fill] (0,0) circle[radius=0.06];
\draw[fill] (1,0) circle[radius=0.06];
\draw[ultra thick] (-1,-2)--(-1,1);
\draw[ultra thick, dashed] (-0,-1)--(-0,1);
\draw[ultra thick] (1,-1)--(1,2);
\draw[ultra thick] (-1,1)--(1,2);
\draw[ultra thick] (-1,-2)--(1,-1);
\draw[fill] (1,-1) circle[radius=0.1];
\draw[fill] (1,2) circle[radius=0.1];
\draw[fill] (-1,-2) circle[radius=0.1];
\draw[fill] (-1,1) circle[radius=0.1];
\draw[fill] (-5,1) circle[radius=0.06];
\draw[fill] (-4,-1) circle[radius=0.06];
\draw[fill] (-1,-1) circle[radius=0.06];
\draw[fill] (-1,0) circle[radius=0.06];
\draw[fill] (1,1) circle[radius=0.06];
\draw[fill] (0,-1) circle[radius=0.06];
\draw[fill] (0,1) circle[radius=0.06];
\node[below] at (-2.5,-2.5) {Figure~4a.};
\node[below] at (-2,-3) {The polygons $P_0^{\gamma,\sigma}$ and $P_1^{\gamma,\sigma}$}; 
\draw[ultra thick](3,0)--(4,0);
\draw[fill] (3,0) circle[radius=0.06];
\draw[fill] (4,0) circle[radius=0.06];
\draw[ultra thick] (3,-2)--(3,2);
\draw[fill] (3,-2) circle[radius=0.1];
\draw[fill] (3,2) circle[radius=0.1];
\draw[ultra thick] (4,-2)--(4,1);
\draw[fill] (4,-2) circle[radius=0.1];
\draw[fill] (4,1) circle[radius=0.1];
\draw[ultra thick] (3,2)--(4,1);
\draw[ultra thick] (3,-2)--(4,-2);
\draw[ultra thick](7,0)--(9,0);
\draw[fill] (7,0) circle[radius=0.06];
\draw[fill] (8,0) circle[radius=0.06];
\draw[fill] (9,0) circle[radius=0.06];
\draw[ultra thick] (7,-1)--(7,1);
\draw[ultra thick, dashed] (8,-2)--(8,1);
\draw[ultra thick] (9,-1)--(9,2);
\draw[ultra thick] (7,1)--(9,2);
\draw[ultra thick] (7,-1)--(8,-2)--(9,-1);
\draw[fill] (9,-1) circle[radius=0.1];
\draw[fill] (9,2) circle[radius=0.1];
\draw[fill] (8,-2) circle[radius=0.1];
\draw[fill] (7,-1) circle[radius=0.1];
\draw[fill] (7,1) circle[radius=0.1];
\draw[fill] (3,1) circle[radius=0.06];
\draw[fill] (3,-1) circle[radius=0.06];
\draw[fill] (4,-1) circle[radius=0.06];
\draw[fill] (8,-1) circle[radius=0.06];
\draw[fill] (8,1) circle[radius=0.06];
\draw[fill] (9,1) circle[radius=0.06];
\node[below] at (6,-2.5) {Figure~4b.};
\node[below] at (6 ,-3) {The polygons $P_0^{\gamma,\delta}$ and $P_1^{\gamma,\delta}$}; 
\end{tikzpicture}
\end{center}

\begin{theor}\label{restheo}
Let $A=(A_0,\ldots, A_n)$ be a tuple of finite sets in $\Z^n$ satisfying the properties given in Definition \ref{resuldef} and $\gamma,\sigma\in(\Z^*)^{|A|}$ be a pair of gradings with strictly positive coordinates and $0$-dimensional support faces $\newton(\resul)_{\gamma}$ and $\newton(\resul)_{\sigma}$. Then the quotient of the coefficients ${r_{\gamma}}$ and ${r_{\sigma}}$ of $\resul(A)$ that are leading with respect to the gradings $\gamma$ and $\sigma$ respectively can be computed as follows:
\begin{equation}\label{formula_res}
\frac{r_{\gamma}}{r_{\sigma}}=(-1)^{\MV(P_0^{\gamma,\sigma},\ldots, P_n^{\gamma,\sigma})}(-1)^{\MV_2(P_0^{\gamma,\sigma},\ldots, P_n^{\gamma,\sigma},\bigl(\begin{smallmatrix}{0}\\{1}\end{smallmatrix}\bigr))}. 
\end{equation}
\end{theor}

\begin{exa}
Using Theorem \ref{restheo}, let us compute the quotient of the coefficients ${r_{\gamma}}$ and ${r_{\sigma}}$  corresponding to the vertices $\bar\gamma$ and $\bar\sigma$  which were considered in \ref{exa1}:

\begin{equation*}
\frac{r_{\gamma}}{r_{\sigma}}=1\cdot (-1)^{\det_2\bigl(\begin{smallmatrix} {1} & {0} & {0}\\{1} & {1} & {1} \end{smallmatrix}\bigr)+\det_2\bigl(\begin{smallmatrix} {0} & {0} & {0}\\{1} & {1} & {1} \end{smallmatrix}\bigr)+\det_2\bigl(\begin{smallmatrix} {1} & {0} & {0}\\{1} & {0} & {1} \end{smallmatrix}\bigr)}=(-1)^0=1.
\end{equation*}

For the coefficients ${r_{\gamma}}$ and ${r_{\delta}}$ corresponding to the vertices $\bar\gamma$ and $\bar\delta$, we obtain

\begin{equation*}
\frac{r_{\gamma}}{r_{\delta}}=1\cdot (-1)^{\det_2\bigl(\begin{smallmatrix} {0} & {1} & {0}\\{1} & {1} & {1} \end{smallmatrix}\bigr)+\det_2\bigl(\begin{smallmatrix} {0} & {0} & {0}\\{1} & {1} & {1} \end{smallmatrix}\bigr)+\det_2\bigl(\begin{smallmatrix} {1} & {0} & {0}\\{0} & {0} & {1} \end{smallmatrix}\bigr)+\det_2\bigl(\begin{smallmatrix} {1} & {1} & {0}\\{0} & {1} & {1} \end{smallmatrix}\bigr)}=(-1)^1=-1.
\end{equation*}

Thus, we obtain the well-known formula for the resultant $\resul=\resul(f,g)$ of the polynomials $f=a_0+a_1x, g= b_0+b_1x+b_2x^2$:
we have $\resul=\pm(a_0^2b_2+a_1^2b_0-a_0a_1b_1)$, just as expected.
\end{exa}

The rest of this Subsection is devoted to the proof of Theorem \ref{restheo}.

\begin{defin}
To the gradings $\gamma, \sigma$, we can associate the {\it  Khovanskii curve} $\mathscr{C}^{\gamma,\sigma}\subset\C^{|A|}$ parametrized by the complex parameter $t\neq 0$ and defined by the following equations: $z_{i,a}=t^{\alpha_{i,a}}+t^{-\beta_{i,a}}$, where $i\in\{0,\ldots,n\}$ and $a\in A_i$.
\end{defin}

Restricting the resultant $\resul$ to the Khovanskii curve $\mathscr{C}^{\gamma,\sigma}$, we obtain a Laurent polynomial in the variable $t$, which we denote by $\phi(t)$. The following statements are obvious.

\begin{utver}
The coefficient of the leading (lowest) term of $\phi(t)$ equals $r_{\gamma}$ ($r_{\sigma}$, respectively).
\end{utver}

\begin{utver}\label{zerophi}
The equality $\phi(t_0)=0$ holds if and only if the point with coordinates $(t_0^{\alpha_{i,a}}+t_0^{-\beta_{i,a}}\mid i\in\{0,\ldots,n\}, a\in A_i)$ belongs to the set $\{\resul=0\}\cap\mathscr C^{\gamma,\sigma}$.
\end{utver}

\begin{rem}\label{res_t}
Note that the polytopes $P_0^{\gamma,\sigma},\ldots, P_n^{\gamma,\sigma}$ are exactly the Newton polytopes of the Laurent polynomials $g_0(x,t),\ldots, g_n(x,t)$, where $$g_i(x,t)=\sum_{a\in A_i}(t^{\alpha_{i,a}}+t^{-\beta_{i,a}})x^a.$$ 
\end{rem}

\begin{proof}[Proof of Theorem \ref{restheo}]
Remark \ref{res_t} implies that the equality (\ref{formula_res}) can be rewritten as follows:
\begin{equation*}
\frac{r_{\gamma}}{r_{\sigma}}=(-1)^{\MV(\newton(g_0),\ldots,\newton(g_n))}(-1)^{\MV_2(\newton(g_0),\ldots,\newton(g_n),\bigl(\begin{smallmatrix}{0}\\{1}\end{smallmatrix}\bigr))}.
\end{equation*}

At the same time, using the classical Vieta's formula, we obtain 
\begin{equation*}
\frac{r_{\sigma}}{r_{\gamma}}=\prod_{\phi(t)=0}t.
\end{equation*}

It follows from Proposition \ref{zerophi} and the Bernstein theorem (see \cite{bernstein} for the details) that
\begin{equation*}
\prod_{\phi(t)=0}t=(-1)^{|\{\resul=0\}\cap\mathscr C^{\gamma,\sigma}|}\prod_{g_0(x,t)=\ldots=g_n(x,t)=0}t=(-1)^{\MV(\newton(g_0),\ldots,\newton(g_n))}\prod_{g_0(x,t)=\ldots=g_n(x,t)=0}t.
\end{equation*}
Then, applying the multivariate Vieta's formula (see Theorem \ref{VietaMain}), we have
\begin{equation*}
(-1)^{\MV(\newton(g_0),\ldots,\newton(g_n))}\prod_{g_0(x,t)=\ldots=g_n(x,t)=0}t=(-1)^{\MV(\newton(g_0),\ldots,\newton(g_n))}(-1)^{\MV_2(\newton(g_0),\ldots,\newton(g_n),\bigl(\begin{smallmatrix}{0}\\{1}\end{smallmatrix}\bigr))},
\end{equation*}
which finishes the proof of the theorem.
\end{proof}

\end{document}